\documentstyle[12pt]{article}
\begin{document}

\

{\large\bf AN ALGEBRA OF SKEW PRIMITIVE ELEMENTS}

\

\centerline{\large\bf V. K. Kharchenko}

\

\

{\it We study the various term operations on the set of
skew primitive elements of Hopf algebras, generated by
skew primitive semi-invariants of an Abelian group of
grouplike elements. All 1-linear binary operations are described and
trilinear and quadrilinear operations are given a detailed
treatment. Necessary and sufficient conditions for the existence of
multilinear operations are specified in terms of the property of
particular noncommutative  polynomials being linearly dependent and of one
arithmetic condition. We dub the conjecture that this condition
implies, in fact, the linear dependence of the polynomials in
question and so is itself sufficient $($ a proof of this conjecture see in "An existence condition for multilinear quantum operations," Journal of Algebra, $217, 1999, 188-228).$} 

\

\centerline{\bf INTRODUCTION}

\

Skew primitive elements in quantum group theory play roughly
the same part as the primitive elements play in the theory of
classical continuous groups. The significance of primitive elements
is determined by the fact that, in any Hopf algebra $H$, the set
$L_1$ of primitive elements form a Lie algebra under the bilinear
term operation $[xy]=xy-yx$, and the subalgebra generated by $L_1$
in $H$ is isomorphic to a universal enveloping (or $p$-enveloping,
if the characteristic $p$ of a ground field is
positive) algebra of the Lie algebra $L_1$.

The study of skew derivations of associative rings, we observe, is
also tightly linked with research on skew primitive elements
of Hopf algebras, since skew primitive elements, indeed,
always act as skew derivations; see [1].

The goal of the present article is to elucidate an algebraic
structure of the set of skew primitive elements of some Hopf
algebras. To take off the ground, it is worth noting that even the
linear structure of that set exhibits itself in a more complex
manner as compared to the way the structure of the set of primitive
elements does. That set forms a {\it comb} (see [2-4]) or, in
other terms, its linear span $L$ is an
Yetter--Drinfeld module (cf. [5, 6]) over a group algebra of the
group $G$ of grouplike elements in a given Hopf algebra. The
multiplicative structure (more precisely, the Lie structure), too,
suffers a sufficient distortion and turns into a partial operation
of variable (quantum) arity; in other words, it splits into the set
of closely related partial operations of different arities. These
operations are the subject matter of our research.

We focus on the Hopf algebras generated by semi-invariants w.r.t. a
commutative group $G$ of grouplike elements acting by conjugations.
We call such Hopf algebras {\it character}. Among them, for
instance, are quantum enveloping Drinfeld--Jimbo algebras;
$G$-universal enveloping algebras of Lie color superalgebras; a quantum
plane; and any Hopf algebra generated by skew primitive
elements, provided that it has exactly $n$ grouplike elements which commute
pairwise, and the ground field contains a primitive $n$th root of
unity.

In quantum group theory spaces of primitive elements of
braided Hopf algebras (cf. [7] and [8]) are currently often treated
as ``quantum'' Lie algebras (see, e.g., [9]). Our present results
apply in research of such ``quantum'' Lie algebras provided that
braiding is defined via the bigrading by a commutative group and by
its character group using the formula
\begin{equation}
(a\underline{\otimes} b)(c\underline{\otimes} d)=
\chi ^c(g_b)^{-1}\cdot ac\underline{\otimes} bd,
\label{braid}
\end{equation}

\noindent because term operations on those ``quantum'' Lie algebras
are defined by the same terms (polynomials) as are the quantum
operations dealt with in the article (see Prop. 4.2 below or
Radford's theorem in [10] on embedding braided Hopf algebras in
ordinary ones via biproduct).

In this article, we give a description of all unary ``quantum''
operations (Thm 5.1) and of all binary operations linear in one of
the variables (Thm. 6.1), and specify a necessary and sufficient
condition for a nonzero $n$-linear operation to exist (Thm. 7.5).
This criterion is then used in Sec. 8 to study in detail trilinear
and quadrilinear operations. Also, we introduce the notion of a
partial main operation of variable arity in terms of which all
quantum operations of degree $\leq 4$ are expressible in the case
of a ground field of characteristic 0.

\

\centerline{\bf 1. BASIC NOTIONS AND EXAMPLES}

\
 
Let $H$ be an arbitrary Hopf algebra with comultiplication $\Delta
$, counity $\varepsilon $, and antipode $S$. Denote by $G$
the set of all grouplike elements
$$ G=\{ g\in H\, | \, \Delta (g)=g\otimes g, \ \ \varepsilon (g)=1\}. $$

\noindent It is well known that $G$ is a group and $S(g)=g^{-1}$,
in which case the linear space generated by $G$ in $H$ is a group
algebra of $G$, that is, distinct grouplike elements are linearly
independent in $H$. For $g\in G$, put
$$ L_g=\{ h\in H\, | \, \Delta (h)=h\otimes 1+g\otimes h\}. $$

\noindent This set forms a linear space over the ground field, and
we call its elements $g$-{\it primitive}, or {\it skew
primitive} if $g$ is not specified. The action of $G$ on $H$ is
defined by conjugations $h^g=g^{-1}hg$. It is easy to see that
linear spaces $L_g$ are independent, that is, their linear span $L$
is the direct sum
$$ L=\sum_{g\in G}\oplus L_g, $$

\noindent in which case the $L$ is invariant under the
above-specified action, and $L_s^g=L_{g^{-1}sg}$. In other words,
$L$ is an Yetter--Drinfeld module (cf. [5, 6]) over a Hopf
subalgebra $k[G]$, which is a group algebra of the group $G$.

\smallskip
{\bf Definition 1.1.} We say that $h\in H$ is a {\it character
element}, or call it a {\it semi-invariant}, if there exists a
character $\chi :G\rightarrow k^*$ such that, for all $g\in G$,
\begin{equation}
g^{-1}hg=\chi (g)h.
\label{act}
\end{equation}

\noindent If $h$ is a nonzero semi-invariant, then the character
$\chi$ is uniquely determined by (2), and we call $\chi$ a {\it
weight} of $h$ and denote it by $\chi ^h$.

\smallskip
{\bf Definition 1.2.} A Hopf algebra $H$ is called {\it character}
if the group $G$ is commutative and $H$ is generated as an
algebra with unity by character skew primitive elements.

The product of two semi-invariants is again a semi-invariant, and
$\chi ^{ab}=\chi ^a \chi ^b$ (if $ab\neq 0$). Therefore,
semi-invariants generate a character Hopf algebra also as a linear
space. Moreover, using (2), we can easily show that nonzero
semi-invariants of different weights are linearly independent. This
means that any character Hopf algebra is graded by the character
group $G^*$ of $G$:
$$ H=\sum _{\chi \in G^*}\oplus H^{\chi }. $$

\smallskip
{\bf Definition 1.3.} Throughout the article, we refer to the
Yetter-Drinfeld module over a group algebra of the Abelian group $G$
as a {\it quantized space} (in view of the fact that in a character
Hopf algebra, that module plays the role of a generating space). In
this way the quantized space is a linear space, graded by an Abelian
group, on which the action of the group is defined in such a way as
to leave homogeneous components invariant.

Consider some examples of character Hopf algebras.

\smallskip
{\bf Example 1.4.} Quantum enveloping algebra 
${\cal {KM}}$. Let $A=||a_{ij}||$ be an arbitrary $n\times n$-matrix,
for which there exist elements $d_1,\ldots, d_n$ such that
$d_ia_{ij}=d_ja_{ji}$ (any Cartan matrix, for instance, has this
property). The Hopf algebra ${\cal {KM}}$ is generated as an
algebra with unity by elements $E_i, F_i, K_i, K_i^{-1}$, $1\leq
i \leq n$, and is defined by the following relations:
$$ K_iK_j=K_jK_i,\ \ \ K_iK_i^{-1}=K_i^{-1}K_i=1, $$
$$ K_i^{-1}E_jK_i=q^{-d_ia_{ij}}E_j, \ \ \ K_i^{-1}F_jK_i=q^{d_ia_{ij}}F_j, $$

\noindent where $q$ is some fixed parameter --- normally, a formal
variable which is freely adjoined to the ground field.
Comultiplication is obtained via the formulas
$$ \Delta (E_i)=E_i\otimes K_i^{-1}+K_i\otimes E_i,\ \ \Delta (F_i)=F_i\otimes K_i^{-1}+K_i\otimes F_i, $$

\noindent by which the counity and antipode are uniquely determined
thus:
$$ \varepsilon(K_i)=1, \ \varepsilon(K_i^{-1})=1,\ \varepsilon(E_i)=\varepsilon(F_i)=0; $$
$$ S(K_i)=K_i^{-1},\ S(E_i)=-q^{-d_ia_{ii}}E_i,\ S(F_i)=-q^{d_ia_{ii}}F_i. $$

In the present example, the group $G$ is generated by elements
$K_i$; its skew primitive generating elements are the
semi-invariants
$$ e_i=E_iK_i, \ f_i=F_iK_i,\ 1-K_i, $$

\noindent whose weights are given by the formulas:
$$ \chi ^{e_i}(K_j)=q^{-d_ja_{ji}}=q^{-d_ia_{ij}},\ \chi ^{f_i}(K_j)=q^{d_ia_{ij}},\ \chi ^{1-K_i}=id. $$

\smallskip
{\bf Example 1.5.} Quantum Drinfeld--Jimbo enveloping algebra
$U_q({\bf g})$. This can be exemplified by a quotient Hopf
algebra of ${\cal{KM}}$ for the case where $A$ is a Cartan matrix
(in particular, $a_{ii}=2$, $a_{ij}\leq 0$; $|a_{ij}|,\, d_i\in
\{ 1,2,3\}$ for $i\neq j$), defined via
\begin{equation}
[E_i,F_j]=\delta _{ij}\left(K_i^2-K_i^{-2}\over q^{2d_i}-q^{-2d_i}
\right)
\label{kom}
\end{equation}

\noindent ($\delta _{ij}$ is a Kronecker symbol), and by the Serre
quantum relations
\begin{equation}
\sum\limits_{\xi =0}^{1-a_{ij}}(-1)^{\xi }
\left(
\matrix{1-a_{ij}\cr \xi
} \right)_{q^{2d_i}}
E_i^{1-a_{ij}-\xi }E_jE_i^{\xi }=0\ \ \ (i\neq j), \label{S1}
\end{equation}
\begin{equation}
\sum\limits_{\xi =0}^{1-a_{ij}} (-1)^{\xi }\left(\matrix{
1-a_{ij}
\cr \xi}
\right)_{q^{2d_i}}
F_i^{1-a_{ij}-\xi }F_jF_i^{\xi }=0\ \ \ (i\neq j), \label{S2}
\end{equation}

\noindent where the parentheses with indices (binomial coefficients)
are given explicitly as values of the polynomials
$$ \left( \matrix{m\cr n} \right) _{t}= {{(t^m-t^{-m})(t^{(m-1)}-t^{-(m-1)})\cdots (t^{(m-n+1)}-t^{-(m-n+1)})}\over {(t-t^{-1})(t^2-t^{-2})\cdots (t^n-t^{-n})}}. $$

\noindent Relations (3) for skew primitive generators take up
the form
\begin{equation}
e_if_j-q^{-2d_ia_{ij}}f_je_i=\delta _{ij}\left(K_i^4-1\over
q^{4d_i}-1\right), \label{komm}
\end{equation}

\noindent whereas Serre $q$-relations (4) and (5) are left fixed
(with $E_i$ and $F_i$ replaced by $e_i$ and $f_i$,
respectively).

\smallskip
{\bf Example 1.6.} Quantum analog for a Lie--Heisenberg algebra.
This is a Hopf subalgebra of $U_q(sl(3))$:
$$ U_q({\cal{H}})=k\langle E_1,E_2,K_1,K_2,K_1^{-1},K_2^{-1} \rangle \subseteq U_q(sl(3)). $$

Since the Cartan matrix of the algebra $sl(3)$ has the form $A=
\left( \matrix{2&-1&0\cr -1&2&-1\cr 0&-1&2}\right) $, the Serre
$q$-relations take up the form
$$ E_i^2E_j+E_jE_i^2=(q^2+q^{-2})E_iE_jE_i $$

\noindent for $i,j\in \{ 1,2\} $, $i\neq j$.

\smallskip
{\bf Example 1.7.} Quantum plane. We have
$$ A^{(2|0)}_q=k\langle g,x \, |\, xg=q gx \rangle $$
$$ G=\langle g\rangle,\ \Delta (x) =x\otimes 1+g\otimes x. $$

\smallskip
{\bf Example 1.8.} Assume that the primitive $n$th root of unity
belongs in the ground field. Then any Hopf algebra with exactly $n$
grouplike elements which commute pairwise will be character provided that
 it is generated by skew primitive elements.

Indeed, $k[G]$, in this case, is a commutative and completely
reducible algebra of dimension $n$, whose irreducible modules all
have dimension 1. In particular, the invariant spaces $L_g$ split
into direct sums of one-dimensional invariant subspaces, which
consist of character elements generating $H$.

\smallskip
{\bf Example 1.9.} Universal $G$-enveloping algebras of Lie color
superalgebras. The concept of a Lie color superalgebra is related to
some fixed Abelian group $G$ and symmetric bicharacter $\lambda
:G\times G\rightarrow k^*$, defined thus:
$$ \lambda (fh,g)=\lambda (f,g)\lambda (h,g), \ \lambda (f,hg)=\lambda (f,h)\lambda (f,g), \ \lambda (f,g)\lambda (g,f)=1. $$

\noindent The linear $G$-graded space $\Lambda =\sum\limits_{g\in
G}\oplus \Lambda _g$ is called a {\it Lie color superalgebra} if it
is augmented with a bilinear operation satisfying the following
properties:
$$ [a,b]=-\lambda (f,g)[b,a],\ \ a\in \Lambda_ f,\ \ b\in \Lambda
_g, \ \ c\in \Lambda _h,$$
$$ \lambda (f,h)[a,[b,c]]+\lambda (h,g)[c,[a,b]]+\lambda (g,h)[b,[c,a]]=0. $$

Any Lie $(G,\lambda )$-color superalgebra $\Lambda $ has a
universal associative enveloping algebra $U$. That is, there
exists a $G$-graded associative algebra $U=\sum \oplus U_g$,
which contains $\Lambda $ as a generating graded subspace $\Lambda
_g\subseteq U_g$, and the operation on homogeneous elements
$a\in \Lambda _f$, $b\in \Lambda _g$ and $\Lambda $ is expressed
via multiplication in $U$ by the formula
$$ [a,b]=ab-\lambda (f,g)ba. $$

\noindent Moreover, $U$ satisfies the categorical universality
condition (for details, see [11-14]). On $U$, we can define the
action of the group $G$ by setting $ a^g=\lambda (f,g)a$, $a\in
\Lambda _f$, and consider a skew group ring $H^{col}=G*U$, which
has the structure of a Hopf algebra with comultiplication defined on
$G$ and $\Lambda $ via $\Delta (g)=g\otimes g$ and $\Delta
(a)=a\otimes 1+f\otimes a$, $a\in \Lambda_f$ (and this is exactly the
Radford biproduct $U(\Lambda )\star {\bf k} [G]$; see [15]).

Now it is easy to see that $H^{col}$ is a character Hopf algebra,
for which
$$ \chi ^a=\lambda (f,\ \, ), \ a\in \Lambda _f;\ \ 
L_g=\Lambda _g\oplus (1-g)k. $$

\newpage

\centerline{\bf 2. QUANTUM VARIABLES AND QUANTUM OPERATIONS}

\

In what follows, we fix an Abelian group $G$ and assume that, in
the Hopf algebras under examination, the $G$ is interpreted by
grouplike elements. In other words, we consider the category of
Hopf algebras $H$ with distinguished homomorphisms $\varphi _H \, :
\,{\bf k} [G]\rightarrow H$.

\smallskip
{\bf Definition 2.1.} A {\it quantum variable} is one to which an
element $g\in G$ and a character $\chi \in G^*$ are associated. In
a Hopf algebra, accordingly, the quantum variable $x=x_g^{\chi }$,
can assume only $g$-primitive semiinvariants of weight $\chi $ 
values only. A character corresponding to $x$ is denoted by 
$\chi ^x$, and a grouplike element --- by $g_x$.

\smallskip
{\bf Definition 2.2.} A {\it quantum operation} in quantum variables
$x_1, \ldots, x_n$ refers to an associative polynomial in $x_1,
\ldots, x_n$ which yields a skew primitive element given any
values of the quantum variables $x_1, \ldots, x_n$ in Hopf
algebras.

In particular, a {\it homogeneous} quantum operation has this form:
$$ [x_1, \ldots, x_n]=\sum\limits_{\pi \in S_n}\alpha _{\pi } x_{\pi (1)}\cdots x_{\pi (n)}, $$

\noindent where $x_1, \ldots, x_n$ are not necessarily distinct
quantum variables. If those variables are mutually distinct (but
not necessarily of different types), the operation is called {\it
multilinear}.

We give some examples of quantum operations.

\smallskip
{\bf 1.} Commutator. If $G$ is a trivial group, then the usual
commutator $xy-yx$ is a quantum operation. If the ground field has
a positive characteristic $p>0$, there exists a nonmultilinear
operation $x^p$. We can show that all other operations (if, of
course, $G={\rm id}$) will be superpositions of these two (this in
essence exhausts the content of the known theorem due to Friedrichs;
see, e.g., [16, Ch. V, Sec. 4]).

\smallskip
{\bf 2.} Skew commutator. Let $x$ and $y$ be quantum variables.
Write $p_{12}=\chi ^x(g_y)$ and $p_{21}=\chi ^y(g_x)$, assuming
that these parameters are related via $p_{12}p_{21}=1$. Then the
skew commutator
\begin{equation}
[x,y]_{p_{12}}=xy-p_{12}yx
\label{skom}
\end{equation}

\noindent is a quantum operation.

In Example 1.5, we have $\chi ^{f_i}(g_{e_j})\chi ^{e_j}(g_{f_i})
=q^{-2d_ia_{ij}}q^{2d_ia_{ij}}=1$, and so the left sides of (6)
are values of the quantum operations and the right ones are
skew primitive constants.

Likewise, the Lie operation in a Lie color superalgebra $\Lambda $
(if $\Lambda $ is assumed embedded in the Hopf algebra $H^{col}$)
is a quantum operation since $ \chi ^a(g_b)\chi ^b(g_a)=\lambda
(f,g)\lambda (g,f)=1$..

It is easy to see that skew commutators essentially exhaust all the
bilinear quantum operations; see Thm 6.1 for $n=1$.

\smallskip
{\bf 3.} Pareigis quantum operation. Let $\zeta $ be a primitive
$n$th root of unity and $x_1, x_2, \dots, x_n$ be quantum
variables such that $\chi ^{x_i}(g_{x_j})\chi ^{x_j}(g_{x_i})=\zeta
^2$. Then
$$ P_n(x_1,\ldots,x_n)= \sum _{\pi \in S_n}(\prod _{i<j\& \pi (i)>\pi (j)} (\zeta ^{-1}\chi ^{x_{\pi (j)}}(g_{x_{\pi (i)}}))\ x_{\pi (1)}\cdots x_{\pi (n)} $$

\noindent is a quantum operation (see [9, Thm. 3.1, p. 147] and the
remarks under Sec. 4 below).

\smallskip
{\bf 4.} Serre quantum operation. Let $x$ and $y$ be such that
$$ \chi ^x(g_y)=q^{2d_ia_{ij}}= \chi ^y(g_x), \ \ \chi ^y(g_y)=q^{4d_i}, $$

\noindent where the parameters $d_i$ and $a_{ij}$ are the same as
in Example 1.5. We can show, then, that the left parts of Serre
quantum relations are values of the following quantum operations:
$$ S_{ij}(x,y)=\sum _{\xi =0}^{1-a_{ij}}(-1)^{\xi } 
\left(\matrix{1-a_{ij}\cr \xi} 
\right)_{q^{2d_i}} y^{1-a_{ij}-\xi }xy^{\xi }$$

\noindent for $x=e_j$ and $y=e_i$ or for $x=f_j$ and $y=f_i$.

These are examples of homogeneous binary quantum operations, linear
in one of the variables. A complete description of such operations
will be given later, under Sec. 6. Now we present the construction
of a ``tensor algebra'' for a quantized space, which makes it possible
to formally operate with polynomials in quantum variables as if with
elements of Hopf algebras.

\

\centerline{\bf 3. FREE ENVELOPING ALGEBRA OF A QUANTIZED SPACE}

\

Let $L=\sum L_g$ be some quantized space. Denote by ${\bf k} \langle
L\rangle $ the tensor algebra of a linear space $L$. If we
distinguish some basis $X$ in $L$, consisting of character
elements, then ${\bf k}\langle L\rangle $ will be a free
associative algebra of $X$ --- in particular, it has a basis
consisting of all words on $X$ (including the empty word equal to
unity). The action of $G$ is uniquely extended to ${\bf k}\langle
L\rangle $, and so we can define the skew group algebra $G* {\bf k}
\langle L\rangle $, on which the structure of a Hopf algebra
arises naturally. We have
$$ \Delta (g)=g\otimes g, \ \varepsilon (g)=1,\ S(g)=g^{-1},\ \ g\in G; $$
$$ \Delta (l)=l\otimes 1+g\otimes l, \ \varepsilon (l)=0,\ S(l)=-g^{-1}l,\ \ l\in L_g.$$

\smallskip
{\bf Definition 3.1.} A {\it free enveloping algebra} of a quantized
space $L$ is the Hopf algebra specified above, and we denote it by
$H\langle L\rangle $, or by $H\langle X\rangle $ if $L$ has some
distinguished basis $X$ consisting of homogeneous character
elements.

Denote the quantized space of skew primitive elements of the
free enveloping algebra $H\langle L\rangle $ by ${\bf
k}_{q-lie}\langle L\rangle$ or by ${\bf k}_{q-lie}\langle X\rangle
$. It is then obvious that $L$ is a quantum subspace in ${\bf
k}_{q-lie}\langle L\rangle $.

Now suppose that some set $X$ of quantum variables is given.
Consider linear spaces $L_x= {\bf k} x$ spanned by the variables
$x$, and denote by $L_g$ the direct sum of all $L_x$ such that
$g_x=g$. Define on $L_g$ the action of $G$ by setting $x^h=\chi
^x(h)x$. In this way the direct sum $L$ of spaces $L_g$ turns
into a quantized space. It is easy to see that quantum operations are
elements of the quantized space $_{q-lie}\langle X\rangle $, in the
representations of which there are no grouplike elements (i.e., they
lie in ${\bf k}\langle X\rangle $). Of course, it is not pointless
to treat all elements ${\bf k}_{q-lie}\langle X\rangle $ as quantum
operations (with constants). This might seem to be even natural
since we have fixed the $G$. Below, however, the reader will see
that this does not in fact give way to any new operations (but for
the constants $1-g$ proper).

A free enveloping algebra is the algebra ${\cal KM}$ from Example
1.4. In this case $G$ is an Abelian group, freely generated by
elements $K_i$, $1\leq i, \leq n$; linear spaces $L_g$ are
either two-dimensional $L_{K_i^2}={\bf k} e_i\oplus {\bf k} f_i$
or zero; the characters are determined from columns of the matrix
$||a_{ij}||$ by setting $\chi ^j (K_i)=q^{d_ia_{ij}}$; and the
action of the group is defined so that $f_j$ turns into a
semi-invariant of weight $\chi ^j$ and $e_j$ turns into a
semi-invariant of weight $(\chi ^j)^{-1}$. Curiously, in order to
construct that quantized space, we need not impose any restrictions on
the matrix $||a_{ij}||$; moreover, if $||a_{ij}||$ is freed of
zero columns, and $q\neq \pm 1$, then, for any character $\chi
\in G^*$ and for every element of $g\in G$, the space $L_g^{\chi
}$ is not more than one-dimensional or, in other words, ${\cal KM}$
does not contain distinct quantum variables of the same type.

We make some trivial but important remarks. First, we bring out the
general form into which a word in quantum variables is expanded
under comultiplication. Let $w=x_1x_2\cdots x_n$; then
$$ \Delta (w)=\Delta (x_1)\Delta (x_2)\cdots \Delta (x_n)=
(x_1\otimes 1+g_{x_1}\otimes x_1)(x_2\otimes 1+g_{x_2}\otimes x_2)
\cdots (x_n\otimes 1+g_{x_n}\otimes x_n).$$

\noindent Removing the parentheses gives
\begin{equation}
\Delta (w)=\sum\limits_{v\in {\cal B}(w)} w|v\otimes v,
\label{w}
\end{equation}

\noindent where ${\cal B}(w)$ denotes the set of all subwords of
$w$ including the empty word (a {\it subword} of $w$ is a word
obtained from $w$ by deleting the letters); $w|v$ is a word
obtained from $w$ by replacing all variables $x_i$ in $v$ by
respective $g_{x_i}$. Put $g_w=g_{x_1}g_{x_2}\cdots g_{x_n}$.
Taking into account that grouplike elements $g_x$ commute pairwise
and that $xg=\chi ^x(g)\cdot gx$, we can reduce (8) to the form
\begin{equation}
\Delta (w)=\sum\limits_{v\in {\cal B}(w)} \alpha _vg_v[w-v]\otimes v,
\label{w1}
\end{equation}

\noindent where $[w-v]$ is a word obtained from $w$ by deleting
$v$, and $\alpha _v$ is a product of all elements of the form
$\chi ^x (g_y)$ for all pairs of variables $x,y$ which occur in
$w$ and are such that $x$ is in $[w-v]$ and $y$ in $v$, and
the quantum variable $x$ occurs in $w$ to the left of $y$.

The above formula can be presented in another form, by writing the
grouplike elements $g_v$ to the right of $[w-v]$ in the left
components of tensors:
\begin{equation}
\Delta (w)=\sum\limits_{v\in {\cal B}(w)} \alpha _v^{\prime
}[w-v]g_v\otimes v, \label{w0}
\end{equation}

\noindent where $\alpha _v^{\prime }$ is a product of all elements
of the form $(\chi ^x (g_y))^{-1}$ for all pairs of variables $y,
x$ which occur in $w$ and are such that $x$ is in $[w-v]$ and
$y$ in $v$, and the quantum variable $x$ occurs in $w$ to the
right of $y$.

Further, on the free enveloping algebra $H\langle X\rangle$ we can
define a degree function $d$ by setting $d(g)=0$, $g\in G$;
$d(x)=1$, $x\in X$. On the tensor product, that function induces
two degree functions:
$$ d_l(w\otimes 1)=d(w),\ d_l(1\otimes w)=0 $$

\noindent and
$$ d_r(w\otimes 1)=0,\ d_r(1\otimes w)=d(w). $$

\noindent Also, we can define the degree $d_+=d_r+d_l$. The tensor
square of a free enveloping algebra has gradings relative to each
one of the degrees. It is worth mentioning that comultiplication
will be homogeneous in view of (9) once we have assumed that
$H\langle X\rangle$ is graded by $d$, and $H\langle X\rangle
\otimes H\langle X\rangle$ --- by $d_+$. In particular, it
follows that all $d$-homogeneous components of quantum operations
are quantum operations themselves, and we can therefore limit our
treatment to $d$-homogeneous quantum operations. It is not hard to
see that the filtration, defined by the degree function $d$,
$$ {\bf k}[G]\subseteq {\bf k}[G]L\subseteq {\bf k}[G]L^2\subseteq \ldots \subseteq {\bf k}[G]L^n\subseteq \ldots $$

\noindent is contained in the coradical filtration; see [17, p. 60].

We will need yet another degree function which is related to some
distinguished variable $x\in X$ and is defined similarly as
follows: $d^{(x)}(x)=1$ and $d^{(x)}(y)=0$, for $y\in X$ and
$y\neq x$, and $d^{(x)}_l(w\otimes v)=d^{(x)}(w)$,
$d^{(x)}_r(w\otimes v)=d^{(x)}(v)$. Clearly, comultiplication will
be homogeneous if we consider the degree $d^{(x)}$ on $H\langle
X\rangle$ and consider $d^{(x)}_+= d^{(x)}_l+d^{(x)}_r$ on
$H\langle X\rangle \otimes H\langle X\rangle$.

Now note that no new grouplike elements arise from a free enveloping
algebra.

\smallskip
{\bf LEMMA 3.2.} {\it Every grouplike element of a free enveloping
algebra of a quantized space belongs to} $G$.

{\bf Proof.} By construction, the basis of a free enveloping algebra
consists of words of the form $gw$, where $g\in G$ and $w$ is a
word in some set $X$ of quantum variables. Then the basis for a
tensor product $H\langle X\rangle \otimes H\langle X\rangle $
consists of tensors of the form $gw\otimes hv$ --- in particular,
those tensors are linearly independent. If $f=\sum \alpha _{gw}gw$
is a grouplike element, then $\Delta f= f\otimes f$; therefore,
\begin{equation}
\sum\limits_{g,w}\alpha _{gw}(g\otimes g)\Delta (w)=
\sum\limits_{g,w}\alpha _{gw}gw\otimes \sum\limits_{g,w}\alpha
_{gw}gw. \label{g}
\end{equation}

\noindent If $w$ is some nonempty word of the greatest length
possible, occurring in the expansion of $f$ with nonzero $\alpha
_{gw}$, then the right-hand side of (11) has the term $\alpha
^2_{gw}gw\otimes gw$, which cannot arise from the left by (9) [or
else in view of the property of comultiplication being
$(d,d_+)$-homogeneous]. Thus, all words in the expansion of $f$ are empty,
that is, $f\in {\bf k} [G]$.

It remains to appeal to the trite fact that all grouplike elements
of a group algebra belong to the initial group. The lemma is proved.

Further, we note that the concept of a quantum operation with
constants brings about nothing new.

\smallskip
{\bf Proposition 3.3.} {\it Every quantum operation with constants of
positive $d$-degree lies in} ${\bf k}\langle X\rangle $.

{\bf Proof.} Let $f=\sum\limits_{h,w}\beta _{hw}hw$ be a skew
primitive element of a free enveloping algebra. Then
$$ f\otimes 1+g_f\otimes f= \sum _{h,w}\beta _{hw}hw\otimes 1+g_f\otimes \sum _{h,w}\beta _{hw}hw= \Delta (\sum _{h,w}\beta _{hw}hw)= $$
$$ \sum _{h,w}\beta _{hw}(h\otimes h)\sum _{v\in {\cal B}(w)} \alpha _vg_v[w-v]\otimes v= \sum _{h,w,v}\beta _{hw} \alpha _vhg_v[w-v]\otimes hv. $$

\noindent Since the linear spaces $1 {\bf k}$, $hL$, $h\in G$,
form a direct sum in the free enveloping algebra, all terms of the
form $...\otimes 1$ should be cancelable, that is,
$$ \sum _{h,w}\beta _{hw}hw\otimes 1=\sum _w\beta _{1w}w\otimes 1. $$

\noindent Because $hw$, $h\in G$, are linearly independent, we
conclude that $\beta _{hw}=0$ for $h\neq 1$. The proposition is
proved.

\

\centerline{\bf 4.  BIGRADED HOPF ALGEBRAS}

\

In quantum group theory, spaces of primitive elements of braided
Hopf algebras are sometimes treated as quantum analogs of Lie
algebras. A {\it braided} Hopf algebra is defined in essentially the
same way as is an ordinary Hopf algebra, the difference being that,
instead of the usual tensor product of algebras, in which the left
and right components commute so that
$$ (1\otimes a)(b\otimes 1)=b\otimes a=(b\otimes 1)(1\otimes a),$$

\noindent we take another product $\cdot $, for which the
commutation rule is given by a Yang--Baxter operator. In this event
the notion of a character Hopf algebra is translated into a concept
of a $G\times G^*$-{\it graded braided} Hopf algebra. Namely, let
${\cal H}$ be an associative algebra graded by the group $G\times
G^*$:
$$ {\cal H}=\sum _{g\in G, \chi \in G^*}\oplus {\cal H}^{\chi }_g. $$

\noindent Redefine multiplication on the tensor product ${\cal
H}\otimes {\cal H}$ of linear spaces by setting
\begin{equation}
(a\otimes b)\cdot (c\otimes d)=(\chi ^c(g_b))^{-1}(ac\otimes bd).
\label{brai}
\end{equation}

\noindent The result is an associative algebra, denoted by ${\cal
H}{\underline{\otimes }}{\cal H}$. Now if, in the definition of a
Hopf algebra, we change the sign $\otimes $ by ${\underline{\otimes
}}$ and assume that the coproduct, counity, and antipode are
homogeneous, we arrive at a definition of the {\it braided bigraded}
Hopf algebra. It will be natural to take the primitive element
$\Delta (a)=a{\underline{\otimes }}1+1{\underline{\otimes }}a$,
lying in the component ${\cal H}^{\chi }_g$, to be the value of
the quantum variable $x=x^{\chi }_g$ in the braided bigraded Hopf
algebra. Accordingly, the definition of a quantum operation
undergoes a slight change, not in content, but in form.

\smallskip
{\bf Definition 4.1.} A {\it braided operation} in quantum variables
$x_1, x_2, \ldots,x_n$ is a polynomial in $x_1, x_2, \ldots,x_n$,
which turns into a primitive element given any values of quantum
variables in braided bigraded Hopf algebras.

\smallskip
{\bf Proposition 4.2.} {\it A homogeneous polynomial $F$ is a 
braided operation if
and only if it is a quantum operation}.

{\bf Proof.} We need only compare the formulas used to compute the
coproducts of monomials in nonbraided [formula (10)] vs. braided
Hopf algebras:
$$ \Delta (w)=\Delta (x_1,x_2,\ldots,x_n)= (x_1\underline{\otimes } 1+1\underline{\otimes } x_1) \cdot (x_2\underline{\otimes } 1+1\underline{\otimes } x_2) \cdot \ \ldots \cdot (x_n\underline{\otimes } 1+1\underline{\otimes } x_n)= $$
$$ \sum _{v\in {\cal B}(w)} \beta _v[w-v]\underline{\otimes } v, $$

\noindent where by
$$ (1\underline{\otimes } x)(y\underline{\otimes } 1)= (\chi
^y(g_x))^{-1}y\underline{\otimes } x,$$

\noindent the coefficients $\beta _v$ are defined in exactly the
same way the $\alpha _v^{\prime }$ are defined in (10), that is,
$\beta _v=\alpha _v^{\prime }$. The proposition is proved.

For a set $X$ of quantum variables, the structure of a bigraded
braided Hopf algebra is naturally defined on the free algebra
$\langle X\rangle $ in a way that quantum variables become
primitive elements. By the above proposition, then, the set of
primitive elements of that algebra is exactly the set of all quantum
operations in $X$.

The above definition of a bigraded braided Hopf algebra is somewhat
different in form from the conventional concept of a
$(G, \lambda )$-{\it graded} Hopf algebra (see, e.g., [15]), in which one grading
${\cal A}=\sum {\cal A}_g$ and a bicharacter $\lambda : G\times
G\rightarrow {\bf k}^*$ are specified and the commutation rule is
defined via
$$ (a\underline{\otimes } b)\cdot (c\underline{\otimes } d)= \lambda (g_b, g_c)(ac\underline{\otimes } bd), $$

\noindent where $b\in {\cal A}_{g_b}$ and $c\in {\cal A}_{g_c}$.
If we put $\chi ^c=\lambda (\ \,, g_c)^{-1}$, the ${\cal A}$
then turns into the bigraded braided Hopf algebra defined above.
Conversely, if, on the group $\underline{G}=G^*\times G$, the
bicharacter is defined via $\lambda (\chi \times g, \chi ^{\prime
}\times g^{\prime }) =\chi ^{\prime }(g)^{-1}$, then the ${\cal
H}$ will turn into a $(\underline{G},\lambda )$-graded Hopf
algebra.

Besides, the Radford biproduct ${\cal A}\star {\bf k} G$ is an
(ordinary) Hopf algebra (see [15, Cor. 3.5]), and if ${\cal A}$ is
generated by primitive elements, then ${\cal A}\star {\bf k} G$ is
a character Hopf algebra whose $g$-primitive elements all have
equal weights. Yet, in an arbitrary (ordinary) Hopf algebra, $g$-primitive
elements do not necessarily all have equal weights ---
this is the reason why we opt for the present definitions.

\

\centerline{\bf 5. UNARY OPERATIONS}

\

{\bf THEOREM 5.1.} {\it For a quantum variable $x$, there exists a
quantum operation $x^n$, $n>1$, if and only if $p= \chi ^x(g_x)$
is a primitive $m$th root of unity, and $n=ml^r$ where $l=1$, if
the characteristic of the ground field ${\bf k}$ is zero, and
$l={\rm char}\,{\bf k}$ if it is positive}.

{\bf Proof.} We have
$$ \Delta (x^n)=(x\otimes 1+g_x\otimes x)^n, $$

\noindent in which case $(x\otimes 1)(g_x\otimes x)=p(g_x\otimes
x)(x\otimes 1)$. Therefore, we can use the so-called quantum
binomial formula, which says that if $XY=pYX$, then
$$ (X+Y)^n=\sum _{k=0}^n\left[ \matrix{n\cr k}\right] _{t=p}Y^kX^{n-k}, $$

\noindent where $\left[ \matrix{n\cr k}\right] _t$ is a polynomial
in $t$, having the following rational representation
$$ \left[ \matrix{n\cr k}\right] _t={t^{[n]}t^{[n-1]}\cdots t^{[n-k+1]}\over t^{[1]}t^{[2]}\cdots t^{[k]}}, $$

\noindent in which by definition, $t^{[s]}=1+t+\cdots +t^{s-1}$,
$t^{[0]}=0$. Using this formula, we obtain
$$ \Delta (x^n)=x^n\otimes 1+g_x^n\otimes x^n+ \sum _{k=1}^{n-1} \left[ \matrix{ n\cr k} \right] _{t=p} g_x^kx^{n-k}\otimes x^k. $$

\noindent In a free enveloping algebra, the elements $1, x, \ldots,
x^n$ are linearly independent; therefore, the $x^n$ will be
primitive if and only if all polynomials $\left[ \matrix{n\cr
k}\right] _t$, $1\leq k\leq n-1$ vanish at $t=p$. In particular,
$p^{[n]}=\left[ \matrix{n\cr 1} \right] _{t=p}=0$, that is, $p^n=1$.

Thus, if $x^n$ is an operation, then $p$ is a primitive $m$th root
of unity for some $m$, and $n$ is divisible by $m$, $n=mq$ (if
$p=1$, we put $m=1$). Consider the coefficient at $\ldots \otimes
x^m$. Note that the following general formula is valid:
\begin{equation}
t^{[ms+d]}=(t^m)^{[s]}t^{[m]}+t^{ms}t^{[d]},
\label{qua}
\end{equation}

\noindent with $0\leq d<m$. Using it yields
$$ \left[ \matrix{ n\cr m} \right] _{t=p}= {t^{[mq]}t^{[mq-1]}\cdots t^{[mq-k+1]}\over t^{[1]}t^{[2]}\cdots t^{[m]}}_{t=p}= {(t^m)^{[q]}t^{[m]}\over t^{[m]}}_{t=p}\cdot {p^{[m-1]}p^{[m-2]}\cdots p^{[1]}\over p^{[1]}p^{[2]}\cdots p^{[m-1]}}=(t^m)^{[q]}|_{t=p}=q. $$

\noindent For the case of characteristic zero, we obtain $\left[
\matrix{ n\cr m} \right] _{t=p}\neq 0$, which is possible only
whenever $m=n$. If $l={\rm char}\,{\bf k}$ is positive, then $q$
is divisible by $l$, and so $n=mlq_1$.

In a similar way, we can compute the coefficient at $\ldots \otimes
x^{ml}$, which in fact equals $q_1$, that is, either $ml=n$,
or $q_1$ is divisible by $l$. Carrying these computations over
and over, we obtain $n=ml^r$. Such routines can be done away with
if we use a simple induction argument along the following lines.

First, $x^m$ is a quantum operation since all coefficients $\left[
\matrix{ m\cr k} \right] _{t=p}$ are such that each has $p^{[m]}=0$
in the numerator and all of their denominators do not vanish. In
addition, $g_{x^m}=g_x^m$ and $\chi ^{x^m}=(\chi^x)^m$; in
particular, $\chi ^{x^m}(g_{x^m})=(p^m)^m=1$. Consider a new
quantum variable $y$ with parameters $\chi ^y=(\chi ^x)^m$ and
$g_y=(g_x)^m$. The (one-dimensional) quantized space generated by
$x^m$ is then isomorphic to a quantized space ${\bf k}y$.
Therefore, free enveloping algebras of these spaces are also
isomorphic under $y\leftrightarrow x^m$. And $y^q$ is a primitive
element since it corresponds to the primitive element $x^n$. Now,
if $m>1$, we can apply the induction hypothesis to $y$, to
eventually obtain $q=l^s$, because $p_y=\chi ^y (y)=1$. But if
$m=1$, that is, $p=1$, then $x^l$ is also a quantum operation,
and a similar inductive step can be taken to treat the quantum
variable $y$ with parameters $\chi ^y=(\chi ^x)^l$ and $g_y=g_x^l$.
We have thus proved the necessity of the condition $n=ml^r$.

Conversely, if $p$ is a primitive $m$th root of unity and if
$n=ml^r$, then $x^n=(\ldots ((x^n)^l)\ldots )^l$, that is, $x^n$,
being a superposition of the quantum operations $x^m$ and $y^l$,
is also a quantum operation. The theorem is proved.

The present theorem shows that there exist two basic types of unary
quantum operations: $x^m$ and $y^l$, where the quantum variables
$x$ and $y$ are such that $\chi ^y(y)=1$ and $\chi ^x(x)$ is a
primitive $m$th root of unity, and all other operations are superpositions of
these two operations. In fact, the two types of operations do not
differ in a crucial respect, since in either case the existence of
an operation $x^n$ follows from the fact that $n$ is a least
number such that $\chi ^x(g_x)^{[n]}=0$. This allows us to
introduce one {\it main unary} (partial) operation $\lbrack\!\lbrack
\ \rbrack\!\rbrack$ on an arbitrary Hopf algebra, defined thus:
$$
\lbrack\!\lbrack a \rbrack\!\rbrack{\buildrel \rm def \over =}
\left\{\begin{array}{ll} a^m, & a \mbox{is a skew primitive
semi-invariant and} \\
& m \mbox{ a minimal number with the property }
\chi ^a(g_a)^{[m]}=0; \\
\mbox{undefined} & \mbox{otherwise.} \end{array}\right.
$$

The main unary operation on an arbitrary braided bigraded Hopf
algebra is defined in exactly the same way: $\lbrack\!\lbrack a
\rbrack\!\rbrack=a^m$ if $a$ is a homogeneous element and $\chi
^a(g_a)^{[m]}=0$, $m$ is minimal.

Clearly, the set of all unary quantum operations form a unary
algebra relative to $\lbrack\!\lbrack \ \, \rbrack\!\rbrack$, and
that algebra is generated by quantum variables.

\newpage

\centerline{\bf 6. BINARY OPERATIONS LINEAR IN ONE OF THE VARIABLES}

\

Let $x$ and $y$ be quantum variables. Denote $p_{12}=\chi ^x(g_y)$,
$p_{21}=\chi ^y(g_x)$, and $p_{22}=\chi ^y(g_y)$.

\smallskip
{\bf THEOREM 6.1.} {\it For quantum variables $x$ and $y$, there
exists a nonzero, quantum, linear in $x$ operation:
\begin{equation}
W(x,y)=\sum\limits_{k =0}^{n}\alpha _ky^{k}xy^{n-k}
\label{2op}
\end{equation}

\noindent if and only if either
\begin{equation}
p_{12}p_{21}=p_{22}^{1-n},
\label{2cond}
\end{equation}

\noindent or $p_{22}$ is a primitive $m>1$th root of unity, $m\mid
n$, and
\begin{equation}
p_{12}^mp_{21}^m=1.
\label{cond}
\end{equation}

If one of these conditions is satisfied, or both, then there exists
a unique, up to multiplication by a scalar, nonzero quantum
operation of degree $n$ w.r.t. $y$, which is linear in $x$.

If condition $(15)$ holds, then the coefficients $\alpha _k$ are
equal to coefficients of the polynomial
\begin{equation}
\prod\limits_{s=0}^{n-1}(t-p_{12}p_{22}^s)=\sum\limits_{k=0}^n
\alpha _kt^k. \label{coef}
\end{equation}

\noindent In this event the operation has the following $q$-commutator
representations:
\begin{equation}
W(x,y)=[\cdots [[xy]_{p_{12}}y]_{p_{12}p_{22}}\cdots
y]_{p_{12}p_{22}^{n-1}}, \label{[1]}
\end{equation}
\begin{equation}
W(x,y)=[\cdots [[xy]_{p_{21}^{-1}}y]_{p_{21}^{-1}p_{22}^{-1}}\cdots
y]_{p_{21}^{-1}p_{22}^{1-n}}.
\label{[2]}
\end{equation}

If $(16)$ holds, then}
\begin{equation}
W(x,y)=[\ldots [[xy^m]_{p_{12}^m}y^m]_{p_{12}^m}\ldots y^m]_{p_{12}^m}.
\label{on}
\end{equation}

{\bf Proof.} First we show that if (15) is satisfied, then the
right-hand sides of equalities (18) and (19) are equal operations
whose coefficients are specified by equality (17).

We consider the sequence of elements $v_0=x$,
$v_{k+1}=[w_ky]_{p_{12}p_{22}^k}$ and use induction to show that
\begin{equation}
d^{(x)}_l(\Delta (v_k)-v_k\otimes 1)=0.
\label{right}
\end{equation}

\noindent Let
$$ \Delta (v_k)-v_k\otimes 1=\sum _{i=0}^k g_xg_y^{n-i}y^i\otimes u_i. $$

\noindent Then
$$ \Delta (v_{k+1})=\Delta (v_k)\Delta (y)-p_{12}p_{22}^k\Delta (y)\Delta (v_k)= $$
$$ (v_k\otimes 1+\sum _{i=0}^k g_xg_y^{n-i}y^i\otimes u_i) (y\otimes 1+g_y\otimes y)- p_{12}p_{22}^k(y\otimes 1+g_y\otimes y) (v_k\otimes 1+\sum _{i=0}^k g_xg_y^{n-i}y^i\otimes u_i). $$

\noindent Removing the parentheses and neglecting terms of the form
$gy^j\otimes \ldots $, $j\geq 0$, (which we can do since their
$d_l^{(x)}$-degrees equal zero), in view of
$v_kg_y=p_{12}p_{22}^kg_yv_k$ we obtain the equality
$$ \Delta (v_{k+1})\equiv v_ky\otimes 1+v_kg_y\otimes y- p_{12}p_{22}^k(yv_k\otimes 1+g_yv_k\otimes y)\equiv (v_ky-p_{12}p_{22}^kyv_k)\otimes 1\equiv v_{k+1}\otimes 1, $$

\noindent as required.

Likewise we consider the sequence of elements $w_0=x$,
$w_{k+1}=[w_ky]_{p_{21}^{-1}p_{22}^{-k}}$ and show that
\begin{equation}
d^{(x)}_r(\Delta (w_k)-g_xg^k_y\otimes w_k)=0.
\label{left}
\end{equation}

\noindent If $$ \Delta (w_k)-g_xg^k_y\otimes w_k=\sum _{i=0}^k
u_i\otimes y^i, $$ then
$$ \Delta (w_{k+1})=\Delta (w_k)\Delta (y)-p_{21}^{-1}p_{22}^{-k}\Delta (y)\Delta (w_k)= $$
$$ (g_xg^k_y\otimes w_k+\sum _{i=0}^k u_i\otimes y^i) (y\otimes 1+g_y\otimes y)- p_{21}^{-1}p_{22}^{-k}(y\otimes 1+g_y\otimes y) (g_xg^k_y\otimes w_k+\sum _{i=0}^k u_i\otimes y^i). $$

\noindent Removing the parentheses and neglecting terms whose
$d_r^{(x)}$-degrees are zero, in view of
$yg_xg_y^k=p_{21}p_{22}^kg_xg_y^ky$ we obtain
$$ \Delta (w_{k+1})\equiv g_xg_y^ky\otimes w_k+g_xg_y^{k+1}\otimes w_ky- p_{21}^{-1}p_{22}^{-k}(yg_xg_y^k\otimes w_k+ g_xg_y^{k+1}\otimes yw_k)\equiv g_xg_y^{k+1}\otimes w_{k+1}, $$

\noindent as required.

Now we show that $v_n=w_n$. To do this, consider an operator
representation of the skew commutators $[wy]_p=w\cdot (R_y-pL_y)$,
where $L_y$ is an operator of left multiplication by $y$ and $R_y$
an operator of right multiplication by $y$. We have
\begin{equation}
v_n=x\cdot (R_y-p_{12}L_y)(R_y-p_{12}p_{22}L_y)\ldots
(R_y-p_{12}p_{22}^{n-1}L_y),
\label{[3]}
\end{equation}
\begin{equation}
w_n=x\cdot (R_y-p_{21}^{-1}L_y)(R_y-p_{21}^{-1}p_{22}^{-1}L_y)\ldots
(R_y-p_{21}^{-1}p_{22}^{1-n}L_y).
\label{[4]}
\end{equation}

\noindent At this point we note that all operators occurring in the
two representations are pairwise commuting, and if condition (15)
is satisfied, then $p_{12}p_{22}^k=p_{21}^{-1}p_{22}^{1-(n-k)}$.
Therefore, the left parts of those two equalities have equal
operators. Thus, $v_n=w_n=W$, and we have
$$ d^{(x)}_+(\Delta (W)-W\otimes 1-g_xg_y^n\otimes W)= $$
$$ d^{(x)}_l(\Delta (v_n) -v_n\otimes 1- g_xg_y^n\otimes v_n)+ d^{(x)}_r(\Delta (w_n)-w_n\otimes 1 -g_xg_y^n\otimes w_n)=0, $$

\noindent which means that $W$ is a skew primitive element.
Now, if we remove the parentheses in (23) and (17) we see that the
coefficients in (14) and in (17) are equal (we can simply replace
$t$ with $R_y\over L_y$, assuming that $R_y$ and $L_y$ are
formal commuting symbols). The sufficiency of condition (15) is
thereby established.

If the second condition of the theorem is satisfied, then $z=y^m$
is a primitive element, and $p_{21}^{\prime }=\chi
^z(g_x)=p_{21}^m$, $p_{12}^{\prime }=\chi ^x(g_z)=p_{12}^m$, and
$p_{22}^{\prime }=\chi^z(g_z)=p_{22}^{m^2}=1$. For $x$ and $z$,
then, condition (15) holds with $n=1$, that is,
$z_1=[xy^m]_{p_{12}^m}$ is a primitive element by the above.
Following up this argument, we will see that the right-hand side of
(20) is a quantum operation.

We argue for the way back. Assume $W(x,y)=\sum\limits_{k
=0}^{n}\alpha _ky^kxy^{n-k}$ is a quantum operation. Then by (8),
\begin{equation}
\sum\limits_{k=0}^n\sum\limits_{\matrix{ v\in {\cal
B}(y^kxy^{n-k})\cr v\neq \emptyset, \ v\neq y^kxy^{n-k}}}\alpha
_ky^kxy^{n-k}|v\otimes v=0. \label{bi1}
\end{equation}

\noindent First we argue for uniqueness. To do this, it suffices to
show that $\alpha _0=0$ implies $W=0$. Suppose, to the contrary,
that $\alpha _0=\ldots =\alpha _{r-1}=0$, $\alpha _r\neq 0$,
$r\geq 1$. Consider all terms of the form $\ldots \otimes
xy^{n-r}$ in (25). The word $xy^{n-r}$ is a subword of just one
word occurring with nonzero coefficient in $W$; therefore, (25)
will have only one term in this form:
$$ \alpha _ry^rg_xg_y^{n-r}\otimes xy^{n-r}.$$

\noindent Consequently, the whole sum cannot vanish. This is a
contradiction, which proves the uniqueness.

Now assume that $\alpha _0\neq 0$ and consider all terms of the
form $\ldots \otimes y$ in (25). The word $y^kxy^{n-k}$ has $n$
different entries of the subword $y$; therefore, the second sum
in (25) has the form
$$ (\sum_{s=0}^k\alpha _ky^sg_yy^{k-s-1}xy^{n-s}+ \sum_{s=0}^{n-k}\alpha _ky^kxy^sg_yy^{n-k-s-1})\otimes y= $$
$$ (\alpha _kp_{22}^{[k]}y^{k-1}xy^{n-k}+\alpha _kp_{22}^{[n-k]}p_{12}p_{22}^k y^kxy^{n-k-1})\otimes y. $$

\noindent Summing all terms of this form over all $k$ and keeping
in mind that different words are linearly independent, we obtain
the system of equalities
\begin{equation}
\alpha _kp_{22}^{[k]}+\alpha _{k-1}p_{22}^{[n-k+1]}p_{12}p_{22}^{k-1}=0,
\ \ k=1, \ldots, n.
\label{bi2}
\end{equation}

\noindent Consider all terms of the form $\ldots \otimes xy^{n-1}$
in (25). The word $xy^{n-1}$ is a subword of just two words:
$yxy^{n-1}$ and $xy^n$. In (25), therefore, there are only two
terms in the desired form: $\alpha _1yg_xg_y^{n-1}\otimes xy^{n-1}$
and $\alpha _0\sum\limits_{s=0}^{n-1}g_xg_y^syg_y^{n-s-1}\otimes
xy^{n-1}$. It follows that $ \alpha _1p_{21}p_{22}^{n-1}+\alpha
_0p_{22}^{[n]}=0$. Comparing this equality with (26) for $k=1$,
$ \alpha _1+\alpha _0p_{22}^{[n]}p_{12}=0$, we obtain
\begin{equation}
\alpha_0p_{22}^{[n]}(-p_{12}p_{21}p_{22}^{n-1}+1)=0. \label{bi3}
\end{equation}

\noindent Therefore, if $p_{22}^{[n]}\neq 0$, then the first
condition of the theorem [see (15)] holds, and $W$ has the desired
form by uniqueness.

To prove the remaining part, we use induction on $n$, that is,
assume that for all lesser values of the parameter $n$, Theorem
6.1 is satisfied in full measure. The basis of induction is the case
$n=1$, which is proved since $p_{22}^{[1]}=1\neq 0$.

Let $m$ be a least number such that $p_{22}^{[m]}=0$, $m>1$.
Then $n$ is divisible by $m$, $n=mq$. By (13), we obtain
$p_{22}^{[ms]}=0$ and $p_{22}^{[ms+d]}\neq 0$ for $1\leq d\leq
m-1$. This allows us to solve the system of equations (26). The
result is $\alpha _{ms+d}=0$ for $1\leq d\leq m-1$, and $\alpha
_{ms}$ are arbitrary parameters, that is, the $W$ takes up the
form
\begin{equation}
W(x,y)=\sum\limits_{s=0}^q\alpha _{ms}y^{ms}xy^{m(q-s)}.
\label{z}
\end{equation}

\noindent By Theorem 5.1, the element $y^m$ is primitive. It is
also obvious that $x$ and $y^m$ generate a free subalgebra in the
free enveloping algebra $H\langle X\rangle $. This means that $x$
and $y^m$, together with the group $G$, generate in $H\langle
X\rangle $ a free enveloping algebra of the quantized space ${\bf
k}x+{\bf k}y^m$. Now consider a new quantum variable $z$ with
parameters $\chi ^z=(\chi^y)^m$ and $g_z=g_y^m$. Then the
quantized spaces ${\bf k}x+{\bf k}y^m$ and ${\bf k}x+{\bf k}z$ will
be isomorphic, and so are their free enveloping algebras under
$x\leftrightarrow x$, $y^m\leftrightarrow z$. Since (28) is a
primitive element in $H\langle {\bf k}x+{\bf k}y^m\rangle $, we
conclude that
$$ N(x,y)=\sum_{s=0}^q\alpha _{ms}z^sxz^{n-s} $$

\noindent is a quantum operation of lesser degree, and so the
theorem applies. We have $\chi ^x(g_z)=p_{12}^m$, $\chi
^z(g_x)=p_{21}^m$, and $\chi ^z(g_z)=(p_{22}^m)^m=1$. In
particular, the latter equality shows that the second condition
cannot be satisfied for $N$, and so the first will hold:
$p_{12}^mp_{21}^m= 1^{1-q}=1$, and
$$ N(x,z)=[\ldots [[xz]_{p_{12}^m}z]_{p_{22}^m}\ldots z]_{p_{22}^m}. $$

\noindent The theorem is complete.

We make some useful remarks. First, it is interesting that the
theorem just proved is --- in form --- not sensitive to the
characteristic of the ground field. This is associated with the fact
that the condition of there being an operation
$[xy^{ml^s}]_{p_{12}^{ml^s}}$ is given by the equality
$p_{12}^{ml^s}p_{21}^{ml^s}=1$, which in the case of
characteristic $l>0$ is equivalent to $p_{12}^{m}p_{21}^{m}=1$ and
ensures the existence of a commutator operation $[xy]_{p_{22}^m}$.

If the variables $x$ and $y$ are such that
$p_{12}=p_{21}=q^{2d_ia_{ij}}$ and $p_{22}=q^{4d_i}$, as in the
Drinfeld--Jimbo algebra, and $q$ is not a root of unity, then the
second condition of the theorem fails. Therefore, the condition of
there being an operation has the form $q^{4d_ia_{ij}}=q^{4d_i(1-n)}$,
whence $n=1-a_{ij}$, that is, only Serre $q$-operations
obtain. If $q$ is a root of unity, and $(q^{4d_i})^{m_i}=1$,
$(q^{4d_ia_{ij}})^{s_i}=1$, $m_i, s_i$ are minimal, then
$\varepsilon =p_{12}^{s_i}=\pm 1$, $\delta = p_{12}^{m_i}=\pm 1$,
and the Serre operations are complemented by $[\ldots
[[xy^{s_i}]_{\varepsilon }y^{s_i}]_{\varepsilon }\ldots
y^{s_i}]_{\varepsilon }$ and by $[\ldots [[S(x,y)y^{m_i}]_{\delta
}y^{m_i}]_{\delta }\ldots y^{m_i}]_{\delta }$. Moreover, by
Theorem 6.1, an operation of degree $n$ relative to $y$ exists
only if $a_{ij}+n-1$ is divisible by $m_i$ or $n$ is divisible
by $s_i$, that is, there are no other (binary 1-linear)
operations by uniqueness.

Furthermore, it is notable that if both conditions are satisfied
together, that is, $p_{12}p_{21}=p_{22}$, $p_{22}$ is a primitive $m$th root of unity, $n=mq$,
the uniqueness yields some identity
\begin{equation}
[xy^m]_{p_{12}^m}=
[\ldots [[xy]_{p_{12}}y]_{p_{12}p_{22}}\ldots y]_{p_{12}p_{22}^{m-1}}.
\label{Óon}
\end{equation}

It is also curious that the quantum operations involved in the
second case are superpositions of the quantum operations in a lesser
degree. Therefore, it seems natural that the specialization of the
multilinear main operation $\lbrack\!\lbrack a,b,\ldots
,b\rbrack\!\rbrack$ should be specified by conditions
$p_{12}p_{21}=p_{22}^{1-n}$ and $p_{22}^{[n]}\neq 0$. Yet, for
the moment we define only a bilinear main operation and observe that
braided bigraded Hopf algebras over a field of characteristic 0, for
which that operation is defined on the whole quantized space of
primitive elements, are universal enveloping algebras of Lie color
superalgebras.

\smallskip
{\bf Definition 6.2.} Main bilinear operation:
$$
\lbrack\!\lbrack a,b \rbrack\!\rbrack{\buildrel \rm def \over =}
\left\{ \begin{array}{ll}
ab-p_{12}ba & \mbox{if } a \mbox{ and } b
\mbox{ are skew primitive character elements and} \\
& p_{12}p_{21}=1,\mbox{ where } p_{12}=\chi ^a(g_b) \mbox{ and }
p_{21}=\chi ^b(g_a); \\
\mbox{undefined} & \mbox{otherwise.} \end{array}\right.
$$

\noindent On braided bigraded Hopf algebras, $\lbrack\!\lbrack a,b
\rbrack\!\rbrack$ is defined in a similar way.

It is easy to see that the above operation satisfies the identities
\begin{equation}
\lbrack\!\lbrack a,b \rbrack\!\rbrack=
-\chi ^a(g_b)\lbrack\!\lbrack b,a \rbrack\!\rbrack,
\label{Y0}
\end{equation}
\begin{equation}
\chi ^a(g_c)\lbrack\!\lbrack a,\lbrack\!\lbrack b,c
\rbrack\!\rbrack \rbrack\!\rbrack
+\chi ^c(g_b)\lbrack\!\lbrack c,\lbrack\!\lbrack a,b
\rbrack\!\rbrack \rbrack\!\rbrack
+\chi ^b(g_a)\lbrack\!\lbrack b,\lbrack\!\lbrack c,a
\rbrack\!\rbrack \rbrack\!\rbrack=0,
\label{Y}
\end{equation}

\noindent subject to the condition that all values $\lbrack\!\lbrack
\,\ \rbrack\!\rbrack $ involved in the representation are
determined.

Now let ${\cal H}$ be a braided bigraded Hopf algebra and assume
that the main operation is defined on all pairs of homogeneous
primitive elements. By linearity, then, $\lbrack\!\lbrack \,\
\rbrack\!\rbrack $ is uniquely determined on the space $\Lambda $
of all primitive elements, and the last two formulas show that
$(\Lambda, \lbrack\!\lbrack \,\ \rbrack\!\rbrack )$ is a Lie
$(\underline{G}, \lambda )$-color superalgebra, where
$\underline{G}$ is a subgroup of $G^*\times G$, generated by
elements of the form $\chi^a\times g_a$, and the bicharacter is
defined by $\lambda (\chi \times g, \chi ^{\prime }\times g^{\prime
})= \chi (g^{\prime })$, while the symmetry of that bicharacter is
implied by the fact that the main operation is total. A standard and
well-known argument will show that in the case of characteristic 0,
${\cal H}$ is itself a universal enveloping algebra of $\Lambda $;
see, e.g., [13].

\newpage

\centerline{\bf 7. MULTILINEAR OPERATIONS}

\

Fix some set of quantum variables $x_1, \ldots,x_n$ and
distinguish one variable $x=x_1$ in it. In what follows, we use
the following notation:
\begin{equation}
\chi ^i=\chi^{x_i};\ g_i=g_{x_i};\ p_{ij}=\chi ^i(g_j);
\ q_k=\prod\limits_{i=1}^{k-1}p_{ik}.
\label{sog}
\end{equation}

\noindent Denote by $S_n$ the permutation group on the set $\{
1,2,\ldots,n\} $, and by $S_n^1$ its subgroup consisting of all
permutations leaving the unity fixed. For our goals, both a
functional and an exponential notation for the action of $S_n$ on
the index set might seem convenient; yet, here we opt for the
second, that is, assume that $i^{(\pi \nu)}=(i^{\pi })^{\nu }=\nu
(\pi (i))$. Write ${\bf \tau }$ to denote the permutation
\begin{equation}
{\bf \tau }=\left(
\matrix{1&2&3&\cdots &n \cr n&n-1&n-2&\cdots &1}\right).
\label{tau}
\end{equation}

\noindent For brevity, we make the convention to write $\pi (A)$ or
$A^{\pi }$ for the permutation $\pi $ and for an arbitrary
expression $A$, meaning that $\pi (A)$ is obtained from $A$ by
applying $\pi $ to each index occurring in $A$ at letters
$p_{ij}$ or at variables $x_i$ --- for instance,
$p_{ij}^{\pi}=p_{\pi (i)\pi (j)}$ or $\pi
(q_k)=\prod\limits_{i=1}^{k-1}p_{\pi (i)\pi (k)}$, but not $\pi
(q_k)=q_{\pi (k)}$. In so doing, we do not require that $S_n$ act
on the ground field. For instance, it might be the case that, in
${\bf k}$, $p_{12}=p_{23}$ is satisfied but
$p_{12}^{(123)}=p_{23}^{(123)}$ is not, that is, this is merely a
notational convention, which is used only unless it leads to
confusion. An arbitrary multilinear polynomial, in accordance with
the above conventions, can be written in the form
\begin{equation}
W(x_1,\ldots,x_n)=\sum\limits_{\pi \in S_n}\alpha _{\pi }\pi (x_1
\cdots x_n). \label{pol}
\end{equation}

\smallskip
{\bf Definition 7.1.} An element $W$ of the free enveloping algebra
$H\langle X\rangle $ of the set $X$ is called {\it left primitive
w.r.t.} $x\in X$ if
$$ d_l^{(x)}(\Delta (W)-W\otimes 1)=0.$$

\noindent Similarly, if we write $H\langle X\rangle _g$ for the
linear span of words $hw$ such that $hg_v=g$, then we say that
$W\in H\langle X\rangle _g$ is {\it right primitive w.r.t.} $x$ if
$$ d_r^{(x)}(\Delta (W)-g\otimes W)=0. $$

\smallskip
{\bf LEMMA 7.2.} {\it If a multilinear polynomial $W$ depending on $x$
is left and right primitive w.r.t. $x$, then the $W$ defines a
quantum operation}.

{\bf Proof.} Since $d_l^{(x)}(g_W\otimes W)=d_r^{(x)}(W\otimes 1)=0$,
we have
$$ d_+^{(x)}(\Delta (W)-W\otimes 1-g_W\otimes W)= $$
$$ d_l^{(x)}(\Delta (W)-W\otimes 1-g_W\otimes W)+ d_r^{(x)}(\Delta
(W)-W\otimes 1-g_W\otimes W)=0.$$

\noindent The lemma is proved.

\smallskip
{\bf THEOREM 7.3.} {\it Polynomial $(34)$ is left primitive w.r.t. $x$ if
and only if it has the following representation}:
\begin{equation}
W=\sum\limits_{\nu \in S^1_n}\beta _{\nu }\nu
([\ldots [[x_1x_2]_{q_2}x_3]_{q_3}\ldots x_n]_{q_n}).
\label{semi}
\end{equation}

{\bf Proof.} First we prove that all summands on the right of (35)
are left primitive w.r.t. $x$. By symmetry, it suffices to
consider only the case $\nu ={\rm id}$. Put $v_1=x_1$, $v_{k+1}=
[v_kx_{k+1}]_{q_{k+1}}$ and use induction to show that all $v_k$
are left primitive w.r.t. $x$. Let $\Delta (v_k)=v_k\otimes
1+\sum u_i\otimes d_i$, where the words $u_i$ are independent of
$x$. Then, if we take into account that $v_kg_{k+1}=\chi
^{v_k}(g_{k+1})g_{k+1}v_k=q_{k+1}g_{k+1}v_k$ and neglect terms
whose $d^{(x)}_l$-degrees equal 0, we obtain
$$ \Delta (v_{k+1})=\Delta (v_k)\Delta (x_{k+1})-q_{k+1} \Delta (x_{k+1})\Delta (v_k)= $$
$$ (v_k\otimes 1+\sum u_i\otimes d_i)(x_{k+1}\otimes 1+g_{k+1}\otimes x_{k+1})- q_{k+1}(x_{k+1}\otimes 1+g_{k+1}\otimes x_{k+1}) (v_k\otimes 1+\sum u_i\otimes d_i)\equiv $$
$$ v_kx_{k+1}\otimes 1+v_kg_{k+1}\otimes x_{k+1}-q_{k+1}
(x_{k+1}v_k\otimes 1+g_{k+1}v_k\otimes x_{k+1})\equiv v_{k+1}\otimes
1.$$

Conversely, assume that $W$ is left primitive w.r.t. $x$.
Consider the element
$$ W^{\prime }=\sum _{\nu \in S_n^1} (-1)^{n-1}\alpha _{{\bf \tau }\nu } \nu ((\prod_{i=2}^n q_i)^{-1} [\ldots [[x_1x_2]_{q_2}x_3]_{q_3}\ldots x_n]_{q_n}). $$

\noindent In the latter formula, we note, each long skew commutator
has exactly one word ending in $x_1=x$, and that word equals
$x_{\nu (n)}\ldots x_{\nu (2)}x_1$ and occurs with coefficient
$(-1)^{n-1}\nu (\prod\limits_{i=2}^nq_i)$. This means that all
words ending in $x$ have equal coefficients in $W$ and in
$W^{\prime }$, that is, the difference $W-W^{\prime }$ has no
words ending in $x$. It remains to show that an element with this
property, left primitive  w.r.t. $x$ equals zero.

Thus, let $\alpha _{\pi }=0$ for $\pi (n)=1$. Denote by $w=uxv$
a word which occurs with nonzero coefficient $\alpha _{\pi}$ in
the representation of $W$ and is such that the subword $v$,
which succeeds $x$, has the least length possible. Then the word
$ux$ is not a subword of the type $[w^{\prime }-v]$
of any word $ w^{\prime }$, except $w$, occurring in the
representation of $W$ with nonzero coefficient. Therefore, the
expansion of $\Delta (W)$ via (8) will show that $\Delta
(W)-W\otimes 1$ has a single tensor of the form $g_vux\otimes v$
with nonzero coefficient. Thus, $W$ cannot be left primitive in
the present case. The theorem is proved.

\smallskip
{\bf THEOREM 7.4.} {\it Polynomial $(34)$ is left primitive w.r.t. $x$ if
and only if it has the representation
\begin{equation}
W=\sum\limits_{\nu \in S^1_n}\beta _{\nu }\nu
([\ldots [[x_1x_2]_{q_2^*}x_3]_{q_3^*}\ldots x_n]_{q_n^*}),
\label{sem}
\end{equation}

\noindent where} $q_k^*=\prod\limits_{i=1}^{k-1}p_{ki}^{-1}$.

The {\bf proof} follows the same line of argument as in the previous
theorem using relations $g_1g_2\ldots g_kx_{k+1}=
q_{k+1}^*x_{k+1}g_1g_2\ldots g_k$ and treating words beginning with
$x$.

\smallskip
{\bf THEOREM 7.5.} {\it For quantum variables $x_1,\ldots,x_n$, there
exists a nonzero quantum multilinear operations if and only if
$\prod\limits_{1\leq i\neq j\leq n}p_{ij}=1$, and the polynomials
\begin{equation}
D_{\nu }{\buildrel \rm def
\over =} \nu ([\ldots [[x_1x_2]_{q_2}x_3]_{q_3}\ldots x_n]_{q_n})-
\nu ([\ldots [[x_1x_2]_{q_2^*}x_3]_{q_3^*}\ldots x_n]_{q_n^*}),
\label{dep}
\end{equation}

\noindent where $q_k^*=\prod\limits_{i=1}^{k-1}p_{ki}^{-1}$ and
$\nu $ runs through $S_n^1$, are linearly dependent in a free
associative algebra. In addition, associated to each linear
dependence $\sum \beta _{\nu }D_{\nu }=0$ is the quantum operation
\begin{equation}
W(x_1,\ldots,x_n)=\sum\limits_{\nu \in S^1_n}\beta _{\nu }\nu
([\ldots [[x_1x_2]_{q_2}x_3]_{q_3}\ldots x_n]_{q_n}).
\label{q}
\end{equation}

\noindent Conversely, every multilinear quantum operation has a
presentation by $(38),$ in which the coefficients $\beta _{\nu }$
determine the linear dependence of} $D_{\nu }$.

{\bf Proof.} Put
\begin{equation}
D_{\nu }^+=
\nu ([\ldots [[x_1x_2]_{q_2}x_3]_{q_3}\ldots x_n]_{q_n}),
\label{q1}
\end{equation}
\begin{equation}
D_{\nu }^-=\nu ([\ldots [[x_1x_2]_{q_2^*}x_3]_{q_3^*}\ldots x_n]_{q_n^*}).
\label{q2}
\end{equation}

\noindent Then $D_{\nu }=D_{\nu }^+-D_{\nu }^-$. Therefore, if
$\sum \beta _{\nu }D_{\nu }=0$, then $\sum \beta _{\nu }D_{\nu
}^+=\sum \beta _{\nu }D_{\nu }^-$. By Theorems 7.3 and 7.4, (38)
is both left and right  primitive w.r.t. $x$, that is, $W$
is a quantum operation by Lemma 7.1. Once we have noted that the
polynomials $D^+_{\nu }$ are linearly independent in a free
algebra, we see that different linear dependences among $D_{\nu }$
correspond to different operations.

Conversely, if $W$ is a quantum operation, we have representation
(38) by Theorem 7.3 and have
\begin{equation}
W=\sum \beta _{\nu }^{\prime }
\nu ([\ldots [[x_1x_2]_{q_2^*}x_3]_{q_3^*}\ldots x_n]_{q_n^*})
\label{q3}
\end{equation}

\noindent by Theorem 7.4. For each permutation $\nu \in S_n^1$,
among all terms on the right of (38), there is only one containing
the word $x_1x_{\nu (2)}\ldots x_{\nu (n)}$ --- this is $\beta
_{\nu }D^+_{\nu }$. In addition, the coefficient at that word
equals $\beta _{\nu }$. By a similar argument, the coefficient at
the same word in (41) equals $\beta _{\nu }^{\prime }$, that is,
$\beta _{\nu }=\beta _{\nu }^{\prime }$, and hence $\sum \beta
_{\nu }D_{\nu }=0$.

We follow the same line to treat words ending in $x$. Comparing
coefficients at $x_{\nu (n)}\dots x_{\nu (2)}x_1$ on the right of
(38) and of (41), we arrive at a system of $(n-1)!$ equalities
$$ \beta _{\nu }(-1)^{n-1}\prod _{k=2}^n(\prod _{i=1}^{k-1}p_{\nu (i)\nu (k)})= \beta _{\nu }(-1)^{n-1}\prod _{k=2}^n(\prod _{i=1}^{k-1}p_{\nu (k)\nu (i)}^{-1}), \ \nu \in S_n^1. $$

\noindent Clearly, all the equalities are equivalent to one:
\begin{equation}
\prod\limits_{1\leq i\neq j\leq n}p_{ij}=1.
\label{X}
\end{equation}

\noindent The theorem is proved.

Our further objective is to prove that condition (42) guarantees
that $D_{\nu }$ are linearly dependent, and hence it is a necessary
and sufficient condition for a nonzero quantum operation to exist.
We set the general solution of this problem aside for a separate
Article (see [18]); here, only the cases $n=3$ and $n=4$ will be discussed in
detail.

\

\centerline{\bf 8. TRILINEAR AND QUADRILINEAR \\ QUANTUM OPERATIONS}

\

{\bf THEOREM 8.1.} {\it For quantum variables $x_1$, $x_2$, and $x_3$,
a nonzero trilinear quantum operation exists if and only if
\begin{equation}
p_{12}p_{21}p_{13}p_{31}p_{23}p_{32}=1.
\label{X3}
\end{equation}

If one of the inequalities
\begin{equation}
p_{12}p_{21}\neq 1, \ p_{13}p_{31}\neq 1, \ p_{23}p_{32}\neq 1
\label{N1}
\end{equation}

\noindent holds, then there exists exactly one (up to multiplication
by a scalar) such operation. If no one of them holds, then all
trilinear operations are linearly expressed in terms of
$\lbrack\!\lbrack x_1,\lbrack\!\lbrack x_2,x_3\rbrack\!\rbrack
\rbrack\!\rbrack$ and $\lbrack\!\lbrack x_2,\lbrack\!\lbrack
x_3,x_1\rbrack\!\rbrack \rbrack\!\rbrack$ via} $(30)$ and $(31).$

{\bf Proof.} For $n=3$, the group $S_3^1$ consists of two
elements, $id$ and (23), for which we have
$$ D_{id}=(p_{21}^{-1}-p_{12})x_2x_1x_3+(p_{31}^{-1}p_{32}^{-1}- p_{13}p_{23})x_3x_1x_2, $$
$$ D_{(23)}=(p_{21}^{-1}p_{23}^{-1}-p_{12}p_{32})x_2x_1x_3+ (p_{31}^{-1}-p_{13})x_3x_1x_2 $$

\noindent [see (37)]. If (43) is met and one of the inequalities
(44) holds (let it be $p_{13}p_{31}\neq 1$ for definiteness), then
$$ (p_{21}^{-1}-p_{12})(p_{31}^{-1}-p_{13})= (p_{31}^{-1}p_{32}^{-1}-p_{13}p_{23}) (p_{21}^{-1}p_{23}^{-1}-p_{12}p_{32}), $$

\noindent and hence
$$ D_{id}-{p_{31}^{-1}p_{32}^{-1}-p_{13}p_{23}\over
p_{31}^{-1}-p_{13}}D_{(23)}=0.$$

\noindent Here, $D_{(23)}\neq 0$, that is, the space generated by
$D_{id}$ and $D_{(23)}$ is one-dimensional, and by Theorem 7.5,
there exists the unique trilinear operation
\begin{equation}
[[x_1x_2]_{p_{12}}x_3]_{p_{13}p_{23}}-
{p_{31}^{-1}p_{32}^{-1}-p_{13}p_{23}\over
p_{31}^{-1}-p_{13}}[[x_1x_3]_{p_{13}}x_2]_{p_{12}p_{32}}.
\label{O3}
\end{equation}

But, if all products $p_{ij}p_{ji}$, $i\neq j$, are equal to
unity, then $D_{id}= D_{(23)}=0$, that is, there exist exactly two
linear dependences between $D_{id}$ and $D_{(23)}$; hence, there
are exactly two linearly independent operations. In this case, on
the other hand, all the three values, $ \lbrack\!\lbrack
x_1,x_2\rbrack\!\rbrack $, $\lbrack\!\lbrack x_1,x_3\rbrack\!\rbrack
$, and $\lbrack\!\lbrack x_2,x_3\rbrack\!\rbrack $, of the main
bilinear operation are defined. Moreover, since $g_{\lbrack\!\lbrack
x_i,x_j\rbrack\!\rbrack } =g_ig_j$ and $\chi ^{\lbrack\!\lbrack
x_i,x_j\rbrack\!\rbrack }= \chi ^i\chi ^j$, we see that $\chi
^{\lbrack\!\lbrack x_i,x_j\rbrack\!\rbrack } (g_k)\chi
^k(g_{\lbrack\!\lbrack x_i,x_j\rbrack\!\rbrack })=1$ holds for
$k\neq i,j$, that is, all possible superpositions, too, are
defined. Among them, by the above argument, only two may be linearly
independent (for example, those specified in the theorem), and the rest are
expressed via them using (30) and (31). The theorem is proved.

Note that if exactly one of the inequalities (44) fails, say,
$p_{12}p_{21}=1$, then the superposition $\lbrack\!\lbrack
x_3,\lbrack\!\lbrack x_1,x_2 \rbrack\!\rbrack \rbrack\!\rbrack$
will be defined; hence, the unique (by Thm. 8.1) quantum operation
will equal that superposition. This circumstance allows us to define
the main trilinear operation in this way:

\smallskip
{\bf Definition 8.2.} Main trilinear operation:
$$
\lbrack\!\lbrack a,b,c\rbrack\!\rbrack{\buildrel \rm def \over =}
\left\{\begin{array}{ll}
[[ab]_{p_{12}}c]_{p_{13}p_{23}}-
{p_{31}^{-1}p_{32}^{-1}-p_{13}p_{23}\over
p_{31}^{-1}-p_{13}}[[ac]_{p_{13}}b]_{p_{12}p_{32}} & \mbox{if }
\prod\limits_{i\neq j}p_{ij}=1\ \mbox{and} p_{ij}p_{ji}\neq 1 \\
& \mbox{for } i\neq j, \mbox{ where } a,b,c \\ &
\mbox{are character skew primitive elements} \\
& \mbox{and }p_{12}=\chi ^a(g_b), \ p_{13}=\chi ^a(g_c), \mbox{ etc.}; \\
\mbox{undefined} & \mbox{otherwise.} \end{array}
\right.
$$

On braided bigraded Hopf algebras, the main trilinear operation is
defined in exactly the same way.

The operation being unique has an implication that if we rename the
variables $x_i\rightarrow x_{\pi (i)}$, then the value of the main
operation on $x_{\pi (1)}$, $x_{\pi (2)}$, $x_{\pi (3)}$ (of
course, it is defined on that sequence since (43) is invariant under
such substitutions) should be linearly expressed via its value on
$x_1$, $x_2$, $x_3$, that is,
\begin{equation}
\lbrack\!\lbrack x_{\pi (1)}, x_{\pi (2)}, x_{\pi (3)} \rbrack\!\rbrack =
\alpha _{\pi }\lbrack\!\lbrack x_1, x_2, x_3 \rbrack\!\rbrack.
\label{Sk}
\end{equation}

\noindent If we compare the coefficients at 
$x_{\pi (1)}x_{\pi (2)}x_{\pi (3)}$ 
on the right- and left-hand sides of (46) we see
that $\alpha _{\pi }=\gamma _{\pi ^{-1}}^{\pi }$ $=$ $\gamma _{\pi }^{-1}$ where $\gamma _{\pi }$ are precisely coefficients in the expansion
$$ \lbrack\!\lbrack x_1, x_2, x_3 \rbrack\!\rbrack= 
\sum \gamma _{\pi } x_{\pi (1)}x_{\pi (2)}x_{\pi (3)};$$

\noindent or, again, by routine computations,
$$ \alpha _{id}=1,\ \alpha _{(123)}={p_{31}-p_{13}^{-1}\over p_{12}-p_{21}^{-1}},\ \alpha _{(132)}= {p_{31}-p_{13}^{-1}\over p_{23}-p_{32}^{-1}},\ 
\alpha _{(13)}=p_{21}p_{32}p_{31}, $$
$$ \alpha _{(12)}=p_{21}p_{23}p_{13} {p_{31}-p_{13}^{-1}\over 
p_{23}-p_{32}^{-1}},\ \ \ \alpha _{(23)}= p_{12}p_{32}p_{13} 
{p_{31}-p_{13}^{-1}\over p_{12}-p_{21}^{-1}}. $$

We pass to the case $n=4$. To the conventions and notation fixed
at the beginning of Sec. 7, we add the following:
\begin{equation}
\{ p_{ij}p_{kl}\ldots p_{rs}\} {\buildrel \rm def \over =}
p_{ij}p_{kl}\cdots p_{rs}-p_{ji}^{-1}p_{lk}^{-1}\ldots p_{sr}^{-1},
\label{-}
\end{equation}

\noindent and for the word $A$ depending on $p_{ij}$, denote by
$\overline{A} $ a word obtained from $A$ by replacing all letters
$p_{ij}$ with $p_{ji}$. These are again merely notational
conventions since we by no means mean that the equality of words in
${\bf k}$ has any bearing on those operators.

\smallskip
{\bf LEMMA 8.3.} {\it Let $C$, $D$, and $E$ be some words in $p_{ij}$.
Then}
\begin{equation}
\{ CE\} \{ D\overline{E}\} -\{ C\} \{ D\} =\{ CD\overline{E}\} \{ E\}.
\label{*}
\end{equation}

{\bf Proof.} Using (47), we rewrite the left- and right-hand sides
of (48) in this way:
$$ \{ CE\} \{ D\overline{E}\} -\{ C\} \{ D\} = (CE-\overline{CE}^{-1})(D\overline{E} -(\overline{D} E)^{-1}- $$
$$ (C-\overline{C}^{-1})(D-\overline{D}^{-1})=CED\overline{E}- \overline{C}^{-1}\overline{E}^{-1}D\overline{E}-CE\overline{D}^{-1}E^{-1}+ $$
$$ \overline{C}^{-1}\overline{E}^{-1}\overline{D}^{-1}E^{-1}-CD+ \overline{C}^{-1}D+C\overline{D}^{-1}-\overline{C}^{-1}\overline{D}^{-1}= $$
$$ CED\overline{E}+\overline{CED}^{-1}E^{-1}-CD-\overline{CD}^{-1}, $$
$$ \{ CD\overline{E}\} \{ E\} =(CD\overline{E}-\overline{CD}^{-1}E^{-1}) (E-\overline{E}^{-1})= $$
$$ CD\overline{E}E-\overline{CD}^{-1}-CD+\overline{CDE}^{-1}E^{-1}. $$

\noindent The lemma is proved.

Note that (43) can be written via braces thus: $\{
p_{12}p_{13}p_{23}\} =0$. Therefore, it might be useful to point
out the following trivial properties of the braces:
\begin{equation}
\{ C\} =0 \rightarrow \{ CD\} =C\{ D\},
\label{1}
\end{equation}
\begin{equation}
\{ C\} =0\ \& \ \{ CD\} =0\ \rightarrow \{ D\} =0.
\label{2}
\end{equation}

\smallskip
{\bf THEOREM 8.4.} {\it For quantum variables $x_1, x_2, x_3, x_4$, a
nonzero quadrilinear quantum operation exists if and only if
\begin{equation}
p_{12}p_{21}p_{13}p_{31}p_{14}p_{41}p_{23}p_{32}p_{24}p_{42}p_{34}p_{43}=1.
\label{X4}
\end{equation}

\noindent If this equality holds, and there is a pair of indices
$i,j$ such that
\begin{equation}
\Gamma _4^{(i,j)}\rightleftharpoons \langle
\{p_{ij}\} \neq 0\ \& \ \{p_{ij}p_{ik}p_{kj}\} \neq 0\
\& \ \{p_{ij}p_{is}p_{sj}\} \neq 0 \rangle,
\label{ga}
\end{equation}

\noindent where $i,j,k,s$ are distinct indices, then there exist
exactly two linearly independent quadrilinear operations.

If condition $\Gamma ^{(ij)}_4$ fails for all $i\neq j$, then all
quadrilinear operations are expressed via the main operation of
ranks $2$ and $3$}.

{\bf Proof.} We seek an element $D_{id}$ in an explicit form.
Expanding the skew commutators in (37) yields
$$
D_{id}=-\{ p_{12}\} x_2x_1x_3x_4-\{ p_{13}p_{23}\} x_3x_1x_2x_4
-\{ p_{14}p_{24}p_{34}\} x_4x_1x_2x_3
$$
\begin{equation}
+\{ p_{12}p_{13}p_{23}\} x_3x_2x_1x_4+\{ p_{12}p_{14}p_{24}p_{34}\} x_4x_2x_1x_3+\{ p_{13}p_{23}p_{14}p_{24}p_{34}\} x_4x_3x_1x_2.
\label{D}
\end{equation}

Now assume that $\beta _{\nu }$ are unknown parameters. Consider
the linear combination $\sum \beta _{\nu } \nu (D_{id})$ and the
coefficients at its distinct words. Setting that combination equal
to zero, we obtain a homogeneous system of twelve equations (equal
to the number of distinct words not beginning with and not ending in
$x_1$) with six unknowns. We show that, under conditions (51) and
$\Gamma _4^{(14)}$, that system has exactly two linearly
independent solutions.

Consider the coefficient at $x_2x_1x_3x_4$. If we apply $\nu \in
S_n^1$, the element $x_1$ will be left fixed; therefore, the word
$x_2x_1x_3x_4$ arises in $\nu (D_{id})$ only from the first three
summands of (53). If it arises from the second, then $\nu (3)=2$,
$\nu (2)=3$, and $\nu (4)=4$, that is, $\nu =(23)$. If it
arises from the third, then $\nu (4)=2$, $\nu (2)=3$, and $\nu
(3)=4$, that is, $\nu =(234)$. Therefore, the whole coefficient
at $x_2x_1x_3x_4$ is equal to
$$ -\{ p_{12}\} \beta _{id} -\{ p_{12}p_{32}\} \beta _{(23)}- \{ p_{12}p_{32}p_{42}\} \beta _{(234)}. $$

\noindent In a similar way, if we compute coefficients at other six
words $\nu (x_2x_1x_3x_4)$, $\nu \in S_n^1$, with $x_1$ holding
second place, we obtain the first group of six equations
\begin{equation}
[-\{ p_{12}\} \beta _{id} -\{ p_{12}p_{32}\} \beta _{(23)}-
\{p_{12}p_{32}p_{42}\} \beta _{(234)}]^{\mu }=0,\ \mu \in S_4^1.
\label{ur}
\end{equation}

\noindent At this point we use the conventions made at the beginning
of Sec. 7, assuming in addition that permutations $\mu $ act on
the indices at $\beta $ by right multiplications: $[\ldots \beta
_{\nu }\ldots ]^{\mu }=\ldots \beta _{\nu \mu }\ldots $.

In exactly the same way we consider the coefficient at
$x_4x_3x_1x_2$. This word arises in $\nu (D_{id})$ from the last
three terms only. If it arises from the last but one term, then
$\nu (4)=4$, $\nu (2)=3$, and $\nu (3)=2$, that is, $\nu
=(23)$. And, if it arises from the fourth, then $\nu (3)=4$,
$\nu (2)=3$, and $\nu (4)=2$, that is, $\nu =(234)$.
Therefore, the coefficient is equal to
$$ \{ p_{13}p_{23}p_{14}p_{24}p_{34}\} \beta _{id} +\{
p_{13}p_{14}p_{34}p_{24}\} \beta _{(23)}+ \{ p_{13}p_{14}p_{34}\}
\beta _{(234)},$$

\noindent and we obtain yet other six equations
\begin{equation}
[\{ p_{13}p_{23}p_{14}p_{24}p_{34}\} \beta _{id}
+\{ p_{13}p_{14}p_{34}p_{24}\} \beta _{(23)}+
\{ p_{13}p_{14}p_{34}\} \beta _{(234)}]^{\mu }=0.
\label{ura}
\end{equation}

\noindent Now compare two equations that correspond to one
permutation $\mu $ in (54) and (55). Using Lemma 8.3 to compute all
the three minors in that system of two equations, we make it sure
that, under condition (51), they all are equal to zero, for example,
$$ \{p_{12}\} \{ p_{13}p_{14}p_{34}\} - \{ p_{12}p_{32}p_{42}\} \{ p_{13}p_{23}p_{14}p_{24}p_{34}\} = - \{ p_{12}p_{13}p_{14}p_{34}p_{23}p_{24}\} \{ p_{32}p_{42}\} =0. $$

\noindent Besides, if some coefficient in (54) equals zero, then by
(50), the corresponding coefficient in (55), too, will be equal to
zero, that is, the whole system of twelve equations is equivalent to
the six in (54).

We order elements of the group $S_4^1$
in this way: $id$,
$(23)$, $(234)$, $(34)$, $(24)$, $(243)$. The matrix of
the system then has the form
$$
\matrix{
\{ p_{12}\} & \{ p_{12}p_{32}\} & \{ p_{12}p_{32}p_{42}\} & 0 & 0 & 0 \cr \ \cr
\{p_{13}p_{23}\} & \{ p_{13}\} & 0 & \{ p_{13}p_{23}p_{43}\} & 0 & 0 \cr \ \cr
0 & 0 & \{ p_{13}\} & 0 & \{ p_{13}p_{43} \} & \{ p_{13}p_{43}p_{23}\} \cr \ \cr
0 & 0 & 0 & \{ p_{12}\} & \{ p_{12}p_{42}p_{32}\} & \{ p_{12}p_{42}\} \cr \ \cr
0 & \{ p_{14}p_{34}p_{24}\} & \{ p_{14}p_{34}\} & 0 & \{ p_{14}\} & 0 \cr \ \cr
\{ p_{14}p_{24}p_{34}\} & 0 & 0 & \{ p_{14}p_{24}\} & 0 & \{ p_{14}\} }
$$

If, in this matrix, we delete the first two columns and the third
and fourth rows, we obtain a triangular submatrix, the leading
diagonal of which has elements $\{ p_{12}p_{32}p_{42}\} $, $\{
p_{13}p_{23}p_{43}\} $, $\{ p_{14}\} $. That is, by condition
$\Gamma ^{(14)}_4$ and remark (50), the corresponding minor is
distinct from zero and the whole system has not more than two
linearly independent solutions. Put $\beta _{id}=1$ and $\beta
_{(23)}=0$, and find one solution for the system of the first two
and last two equations:
$$ \beta _{id}=1,\ \beta _{(23)}=0,\ \beta _{(234)}=-{\{ p_{12}\} \over \{p_{12}p_{32}p_{42}\} },\ \beta _{(34)}=-{\{ p_{13}p_{23}\} \over \{ p_{13}p_{23}p_{43}\} }, $$
\begin{equation}
\beta _{(24)}={\{ p_{12}\} \{ p_{14}p_{34}\} \over
\{ p_{14}\} \{ p_{12}p_{32}p_{42}\} },\
\beta _{(243)}=-{\{ p_{43}\} \{ p_{14}p_{24}p_{34}p_{13}p_{23}\} \over
\{p_{14}\} \{ p_{13}p_{23}p_{43}\} }.
\label{be4}
\end{equation}

\noindent Using Lemma 8.3, we verify whether these values are
solutions for the third equation:
$$ -\{ p_{13}\} {\{ p_{12}\} \over \{p_{12}p_{32}p_{42}\} }+ \{ p_{13}p_{43}\} {\{ p_{12}\} \{ p_{14}p_{34}\} \over \{ p_{14}\} \{ p_{12}p_{32}p_{42}\} }- $$
$$ \{ p_{13}p_{43}p_{23}\} {\{ p_{43}\} \{ p_{14}p_{24}p_{34}p_{13}p_{23}\} \over \{p_{14}\} \{ p_{13}p_{23}p_{43}\} }= -{\{ p_{12}\} \over \{p_{12}p_{32}p_{42}\} \{ p_{14}\}} (\{ p_{13}\} \{ p_{14}\} - $$
$$ \{ p_{13}p_{43}\} \{ p_{14}p_{34}\} )- {\{ p_{43}\} \{ p_{14}p_{24}p_{34}p_{13}p_{23}\} \over \{p_{14}\} }= $$
$$ {\{ p_{43}\} \over \{p_{12}p_{32}p_{42}\} \{ p_{14}\} } (\{ p_{12}\}\{ p_{13}p_{14}p_{34}\} - \{p_{12}p_{32}p_{42}\} \{ p_{14}p_{24}p_{34}p_{13}p_{23}\} )= $$
$$ -{\{ p_{43}\} \{ p_{12}p_{13}p_{14}p_{34}p_{23}p_{24}\} \{ p_{32}p_{42}\} \over \{p_{12}p_{32}p_{42}\} \{ p_{14}\} }=0. $$

\noindent Likewise for the fourth equation (with the ``$-$'' sign):
$$ \{ p_{12}\} {\{ p_{13}p_{23}\} \over \{ p_{13}p_{23}p_{43}\} }- \{p_{12}p_{32}p_{42}\}{\{p_{12}\} \{ p_{14}p_{34}\} \over \{p_{12}p_{32}p_{42}\} \{ p_{14}\} }+ $$
$$ \{ p_{12}p_{42}\}{\{ p_{43}\} \{ p_{14}p_{24}p_{34}p_{13}p_{23}\} \over \{ p_{14}\} \{p_{13}p_{23}p_{43}\} }= { \{ p_{12}\} \over \{ p_{14}\} \{ p_{13}p_{23}p_{43}\} } (\{ p_{13}p_{23}\} \{ p_{14}\} - $$
$$ \{ p_{13}p_{23}p_{43}\}\{ p_{14}p_{34}\} ) +\ldots = -{ \{ p_{12}\} \{ p_{13}p_{23}p_{14}p_{34}\} \{ p_{43}\} \over \{ p_{14}\} \{ p_{13}p_{23}p_{43}\} }+ $$
$$ \{ p_{12}p_{42}\}{\{ p_{43}\} \{ p_{14}p_{24}p_{34}p_{13}p_{23}\} \over \{ p_{14}\} \{p_{13}p_{23}p_{43}\} }= { \{ p_{43}\} \over \{ p_{14}\} \{ p_{13}p_{23}p_{43}\} } (-\{ p_{12}\} \{ p_{13}p_{23}p_{14}p_{34}\} + $$
$$ \{ p_{12}p_{42}\} \{ p_{14}p_{24}p_{34}p_{13}p_{23}\} )= ={\{ p_{43}\} \{ p_{12}p_{13}p_{23}p_{14}p_{34}p_{24}\} \{ p_{42}\} \over \{ p_{14}\} \{ p_{13}p_{23}p_{43}\} }=0. $$

\noindent Thus, by Theorem 7.5, the computed values of $\beta _{\nu
}$ determine the quadrilinear operation
\begin{equation}
\lbrack\!\lbrack x_1, x_2, x_3, x_4 \rbrack\!\rbrack =
\sum \beta _{\nu }D^+_{\nu }.
\label{G4}
\end{equation}

\noindent Since $\beta _{id}=1$, $\beta _{(23)}=0$, and the word
$\nu (x_1x_2x_3x_4)$ in (57) occurs only in the summand $D^+_{\nu
}$, we see that the coefficient at $x_1x_2x_3x_4$ in the
expansion (34) of the polynomial $\lbrack\!\lbrack x_1, x_2, x_3,
x_4 \rbrack\!\rbrack $ equals 1, and the coefficient at
$x_1x_3x_2x_4$ is zero.

Consider a sequence of quantum variables $y_1=x_1$, $y_2=x_3$,
$y_3=x_2$, and $y_4=x_4$. This sequence satisfies both conditions
(51) and $\Gamma _4^{(14)}$; hence, by the above, there exists a
quantum operation $\lbrack\!\lbrack y_1, y_2, y_3, y_4
\rbrack\!\rbrack = \lbrack\!\lbrack x_1, x_3, x_2, x_4
\rbrack\!\rbrack $ such that the coefficient at $x_1x_3x_2x_4$
equals 1 and the one at $x_1x_2x_3x_4$ equals 0. In this way
$\lbrack\!\lbrack x_1, x_3, x_2, x_4 \rbrack\!\rbrack $ supplies
the second solution for the system under consideration, which proves
the first part of the theorem.

Now assume that condition $\Gamma ^{(ij)}_4$ is not satisfied for
any pair $i\neq j$. We call the set of quantum variables $Y$
{\it conforming} if condition (42) is satisfied for it. The failure
of condition $\Gamma _4^{(ij)}$ will mean, then, that the variables
$x_i$ and $x_j$ enter some two- or three-element conforming subset.
If the pair $x_i, x_j$ is itself conforming, then the value
$\lbrack\!\lbrack x_i, x_j\rbrack\!\rbrack $ is defined, and the
set $\lbrack\!\lbrack x_i, x_j\rbrack\!\rbrack $, $x_k$, $x_l$
too is conforming. Therefore, one of the superpositions
$\lbrack\!\lbrack \lbrack\!\lbrack x_i, x_j\rbrack\!\rbrack,
x_k,x_l\rbrack\!\rbrack $ or $\lbrack\!\lbrack \lbrack\!\lbrack
\lbrack\!\lbrack x_i, x_j\rbrack\!\rbrack, x_k\rbrack\!\rbrack,
x_l\rbrack\!\rbrack $ is determined. Similarly, if the triple $x_i,
x_k, x_j$ is conforming, then either $\lbrack\!\lbrack
\lbrack\!\lbrack x_i, x_j, x_k\rbrack\!\rbrack, x_l\rbrack\!\rbrack
$ or $\lbrack\!\lbrack \lbrack\!\lbrack \lbrack\!\lbrack x_i,
x_j\rbrack\!\rbrack, x_k\rbrack\!\rbrack, x_l\rbrack\!\rbrack $
is defined.

We turn on to consider the possible cases where the six conditions
$\Gamma ^{(ij)}_4$ are all fallible.

\smallskip
{\bf 1.} All two-element subsets are conforming. The system (54) has
only zero coefficients, and by Theorem 7.5, we then find six
linearly independent quantum operations $D_{\nu }^+$, $\nu \in
S^1_4$,
$$ D_{\nu }^+= \nu (\lbrack\!\lbrack \lbrack\!\lbrack \lbrack\!\lbrack x_1, x_2\rbrack\!\rbrack, x_3\rbrack\!\rbrack, x_4\rbrack\!\rbrack ). $$

\smallskip
{\bf 2.} All four three-element subsets are conforming. In view of
the above, we can assume that one of the two-element subsets is not
conforming. Suppose $\{ p_{12}\} \neq 0$. Then the system (54)
splits into three pairs of equations: $id$, $(23)$; $(243)$,
$(24)$; $(34)$, $(243)$. Here, the first and third pairs have
rank 1 and the second has rank $\leq 1$. Thus, if at least one of
the inequalities $\{ p_{14}\} \neq 0$, $\{ p_{13}\} \neq 0$, or
$ \{ p_{43}\} \neq 0$ holds, then the whole system has exactly
three solutions, and these, in accordance with Theorem 7.5, yield
the following three operations:
$$ \lbrack\!\lbrack \lbrack\!\lbrack x_1, x_2, x_3\rbrack\!\rbrack
^{\prime }, x_4\rbrack\!\rbrack;\ \lbrack\!\lbrack \lbrack\!\lbrack
x_1, x_3, x_4\rbrack\!\rbrack ^{\prime }, x_2\rbrack\!\rbrack;\
\lbrack\!\lbrack \lbrack\!\lbrack x_1, x_2, x_4\rbrack\!\rbrack
^{\prime }, x_3\rbrack\!\rbrack. $$

\noindent Here, $\lbrack\!\lbrack \rbrack\!\rbrack ^{\prime }$
denotes the ternary operation whose uniqueness is asserted by
Theorem 8.1, that is, it is either the main operation or a
superposition of the form $\lbrack\!\lbrack \lbrack\!\lbrack \,\
\rbrack\!\rbrack,\ \rbrack\!\rbrack $. There then exists one more
superposition $\lbrack\!\lbrack \lbrack\!\lbrack x_2, x_3,
x_4\rbrack\!\rbrack, x_1\rbrack\!\rbrack $, which should be
linearly expressed in terms of the solutions that we have found.
Consequently, using the fact that the bilinear (30) and trilinear
(46) operations are symmetric, we arrive at an analog of the Jacobi
identity
\begin{equation}
\sum\limits_{k=0}^3\xi _k\sigma ^k(\lbrack\!\lbrack
\lbrack\!\lbrack x_1, x_2, x_3\rbrack\!\rbrack , x_4\rbrack\!\rbrack
)=0, \label{Y4}
\end{equation}

\noindent where $\sigma =(1234)$ is a cyclic permutation, the
coefficients $\xi _k$ are uniquely determined up to multiplication
by a common scalar, and all values of the main operation are assumed
determined.

If $\{ p_{14}\} = \{ p_{13}\} =\{ p_{43}\} =0$, then the second
pair of equations disappears, and instead of $\lbrack\!\lbrack
\lbrack\!\lbrack x_1, x_3, x_4\rbrack\!\rbrack ^{\prime },
x_2\rbrack\!\rbrack $, there appear two operations: $D^+_{(234)}=
\lbrack\!\lbrack \lbrack\!\lbrack \lbrack\!\lbrack x_1, x_3
\rbrack\!\rbrack x_4, \rbrack\!\rbrack x_2\rbrack\!\rbrack $ and
$D^+_{(24)}= \lbrack\!\lbrack \lbrack\!\lbrack \lbrack\!\lbrack x_1,
x_4 \rbrack\!\rbrack x_3, \rbrack\!\rbrack x_2\rbrack\!\rbrack $.

\smallskip
{\bf 3.} Three three-element subsets are conforming. Condition (51)
then implies that the fourth subset is also conforming.

\smallskip
{\bf 4.} Exactly two three-element subsets are conforming. To be
specific, let $\{ p_{12}p_{14}p_{24}\} =$ $\{ p_{13}p_{14}p_{34}\}
=0$, $\{ p_{12}p_{23}p_{13}\} \neq 0$, $\{
p_{23}p_{24}p_{34}\}\neq 0$. Since condition $\Gamma _4^{(23)}$
fails and the two triples involved are not conforming, we have $\{
p_{23}\} =0$. If we write (51) in the form $\{
p_{12}p_{14}p_{24}p_{13}p_{23}p_{43}\} =0$, by formulas (49) and
(50), we obtain $0=\{ p_{13}p_{23}p_{43}\} =p_{23}\{ p_{13}p_{43}\}
$, and similarly $0=\{ p_{12}p_{32}p_{42}\} =p_{32}\{
p_{12}p_{42}\} $. In other words, $\{ p_{13}p_{43}\} =\{
p_{12}p_{42}\} =\{ p_{23}\} =0$, and again condition (51) yields
$\{ p_{14}\} =0$. In the matrix of (54), in particular, the last
two columns will disappear, and the minor corresponding to the first
four rows and columns will equal $\{ p_{12}p_{23}p_{13}\} \{
p_{13}\} \{ p_{13}\} \{ p_{12}\}$.

Now if $\{ p_{13}\} \{ p_{12}\} \neq 0$, then the whole system has
rank 4 and its solutions are determined by arbitrary values of
$\beta _{(24)}$ and $\beta _{(243)}$, that is, we obtain two
operations, $D^+_{(24)}$ and $D^+_{(243)}$,
\begin{equation}
\lbrack\!\lbrack \lbrack\!\lbrack \lbrack\!\lbrack x_1, x_4
\rbrack\!\rbrack x_3,
\rbrack\!\rbrack x_2\rbrack\!\rbrack; \ \
\lbrack\!\lbrack \lbrack\!\lbrack \lbrack\!\lbrack x_1, x_4
\rbrack\!\rbrack x_2,
\rbrack\!\rbrack x_3\rbrack\!\rbrack,
\label{B2}
\end{equation}

\noindent in terms of which all other operations defined in the
present case are expressible:
$$ \lbrack\!\lbrack \lbrack\!\lbrack x_1, x_2, x_4 \rbrack\!\rbrack, x_3 \rbrack\!\rbrack,\ \lbrack\!\lbrack \lbrack\!\lbrack x_1, x_3, x_4 \rbrack\!\rbrack, x_2 \rbrack\!\rbrack,\ \lbrack\!\lbrack x_1,\lbrack\!\lbrack x_2, x_3 \rbrack\!\rbrack, x_4 \rbrack\!\rbrack,\ \lbrack\!\lbrack \lbrack\!\lbrack x_1, x_4 \rbrack\!\rbrack, x_2, x_3 \rbrack\!\rbrack. $$

If $\{ p_{13}\} \{ p_{12}\} =0$, in view of the initial conditions
being symmetric under the permutation $2\leftrightarrow 3$, it
suffices to consider the case $\{ p_{13}\} =0$. We have $\{
p_{12}p_{24}\} =\{ p_{13}\} =\{ p_{14}\} =\{ p_{23}\} =\{ p_{34}\}
=0$, $\{ p_{12}\} \neq 0$, $\{ p_{24}\} \neq 0$ (if not all
pairs are conforming). And we face only one additional solution
$D^+_{(234)}= \lbrack\!\lbrack \lbrack\!\lbrack \lbrack\!\lbrack
x_1, x_3 \rbrack\!\rbrack x_4, \rbrack\!\rbrack x_2\rbrack\!\rbrack
$ since the minor corresponding to the first, fourth, and sixth
rows and to the first, second, and fourth columns is not equal to
zero.

\smallskip
{\bf 5.} Only one three-element subset is conforming. Let it be
$x_2, x_3, x_4$. Then the failure of conditions $\Gamma ^{(12)}$,
$\Gamma ^{(13)}$, and $\Gamma ^{(14)}$ implies that $\{
p_{12}\} =\{ p_{13}\} =\{ p_{14}\} =0$. In this case $\{ p_{34}\}
\neq 0$, since otherwise the triple $x_1, x_3, x_4$ would be
conforming. Similarly, $\{ p_{24}\} \neq 0$ and $\{ p_{23}\} \neq
0$. These imply $\{ p_{23}p_{24}\} \neq 0$, $\{ p_{23}p_{34}\}
\neq 0$, and $\{ p_{24}p_{34}\} \neq 0$ since, for instance,
condition $\{ p_{23}p_{24}\} = 0$, combined with $\{ p_{1i}\}=0$,
$i=2, 3, 4$, and (51), yields $\{ p_{34}\} =0$. Under these
conditions, the system splits into three pairs of rank 1 equations:
$(23)$, $(243)$; $id$, $(24)$; $(234)$, $(34)$. The
first pair agrees with the operation
$$ D^+_{id}-{\{ p_{23}p_{43}\} \over \{ p_{23}\} }D^+_{(34)}= \lbrack\!\lbrack \lbrack\!\lbrack x_1, x_2 \rbrack\!\rbrack, x_3, x_4 \rbrack\!\rbrack, $$

\noindent and the two other operations result from substitutions
(23) and (24). All other superpositions defined in the present case
are linearly expressed via these three. Specifically, we have an
identity of the form
$$ \lbrack\!\lbrack x_1,\lbrack\!\lbrack x_2, x_3, x_4 \rbrack\!\rbrack \rbrack\!\rbrack = \xi _1\lbrack\!\lbrack \lbrack\!\lbrack x_1, x_2 \rbrack\!\rbrack, x_3, x_4 \rbrack\!\rbrack + \xi _2\lbrack\!\lbrack x_2,\lbrack\!\lbrack x_1, x_3 \rbrack\!\rbrack, x_4 \rbrack\!\rbrack + \xi _3\lbrack\!\lbrack x_2, x_3, \lbrack\!\lbrack x_1, x_4 \rbrack\!\rbrack \rbrack\!\rbrack. $$

\smallskip
{\bf 6.} No one of the three-element subsets is conforming. Then
two-element subsets cannot all be conforming; therefore, one of the
conditions $\Gamma ^{(ij)}_4$ is satisfied.

The theorem is proved.

\smallskip
{\bf Definition 8.5.} Under condition $\Gamma _4^{(14)}$, the {\it
main quadrilinear operation} is defined by
$$ \lbrack\!\lbrack a_1,a_2,a_3,a_4 \rbrack\!\rbrack = \sum _{\nu (1)=1} \beta _{\nu }[[[a_1a_{\nu (2)}]_{\nu (q_2)}a_{\nu (3)}]_{\nu (q_3)}a_{\nu (4)}]_{\nu (q_4)}, $$

\noindent where $a_i$ are skew primitive character elements,
$p_{ij}=\chi ^{a_i}(g_{a_j})$, $\nu
(q_i)=\prod\limits_{k=1}^{i-1}p_{\nu (k) \nu (i)}$, and the
coefficients $\beta _{\nu }$ are given as in (56).

If no proper subset of the set $x_1, x_2, x_3, x_4$ is conforming,
then all conditions $\Gamma ^{(ij)}_4$ are satisfied. Therefore,
all possible $4!=24$ permutation variants $\lbrack\!\lbrack x_{\pi
(1)},x_{\pi (2)},x_{\pi (3)},x_{\pi (4)} \rbrack\!\rbrack $, $\pi
\in S_4$, are determined, and by Theorem 8.4, they all are
expressible via any pair of them. In order to find that
representation, we write the main operation in the form
\begin{equation}
\lbrack\!\lbrack x_1,x_2,x_3,x_4
\rbrack\!\rbrack=
\sum \alpha _{\pi }x_{\pi (1)}x_{\pi (2)}x_{\pi (3)}x_{\pi (4)},
\label{po}
\end{equation}

\noindent where $\alpha _{\pi }$ are particular rational functions
in $p_{ij}$, obtained by expanding the skew commutators in
Definition 8.5. We have already mentioned that $\alpha _{id}=1$ and
$\alpha _{(23)}=0$. Given an arbitrary replacement $x_i\rightarrow
x_{\mu (i)}$, $\mu \in S_4$, we obtain
$$ \lbrack\!\lbrack x_{\mu (1)},x_{\mu (2)},x_{\mu (3)},x_{\mu (4)} \rbrack\!\rbrack = \sum (\alpha _{\pi })^{\mu } x_{\mu (\pi (1))}x_{\mu (\pi (2))}x_{\mu (\pi (3))}x_{\mu (\pi (4))}. $$

\noindent On the right-hand side of the latter equality, the
coefficient at $x_1x_2x_3x_4$ equals $\alpha ^{\mu }_{\mu ^{-1}}$
and the one at $x_1x_3x_2x_4$ equals 
$\alpha ^{\mu }_{(23)\mu ^{-1}}$. 
Therefore, we have a formula that replaces the
twisted symmetry in (46):
\begin{equation}
\lbrack\!\lbrack x_{\mu (1)},x_{\mu (2)},x_{\mu (3)},x_{\mu (4)}
\rbrack\!\rbrack =
\alpha ^{\mu }_{\mu ^{-1}}\lbrack\!\lbrack x_1,x_2,x_3,x_4
\rbrack\!\rbrack +
\alpha ^{\mu }_{(23)\mu ^{-1}}\lbrack\!\lbrack x_1,x_3,x_2,x_4
\rbrack\!\rbrack.
\label{Sim}
\end{equation}

Clearly, the same trick will help us find an expression for any
quadrilinear operation in terms of the main operation since
coefficients in the expansion are equal to those at $x_1x_2x_3x_4$
and $x_1x_3x_2x_4$. For the Pareigis operation $P_4$, for
instance, we have
\begin{equation}
P_4=\lbrack\!\lbrack x_1,x_2,x_3,x_4
\rbrack\!\rbrack +\zeta ^{-1}p_{23}\lbrack\!\lbrack x_1,x_3,x_2,x_4
\rbrack\!\rbrack.
\label{Par}
\end{equation}

\noindent By definition, then, $p_{ij}p_{ji}=-1$; therefore, $\{
A\} =2A$ for words of odd length and $\{ A\} =0$ for words of even
length. That is, the main operation has the following
representation:
$$ \lbrack\!\lbrack x_1,x_2,x_3,x_4 \rbrack\!\rbrack =D^+_{id}-p_{23}p_{24}D^+_{(234)}-p_{24}p_{34}D^+_{(243)}. $$

In conclusion we note that any (not multilinear) operation admits a
full and partial linearizations. On identifying variables in the
linearized operation, we obtain the initial operation multiplied by
an integer dividing $n!$. Therefore, if the ground field has
characteristic zero, all operations will be expressed via multilinear
ones. If the characteristic of ${\bf k}$ is distinct from 2, 3,
then all operations of degree $\leq 4$ are expressed via the main
operation of variable arity, defined in the article. The picture
changes if the characteristic equals 2 or 3. Assume, for instance,
that it equals 2. Then $\lbrack\!\lbrack x, y, y \rbrack\!\rbrack
=0$ if the left-hand side is determined, and by Theorem 6.1, there
still exists a nonzero operation
$xy^2+(p_{12}+p_{12}p_{22})yxy+y^2x$.

\

Acknowledgement. I am extremely indebted to the
participants of Shirshov Seminar on Ring Theory 
of the Institute of Mathematics RAS held in July, 1997, particulary
L.A. Bokut', I. P. Shestakov, V. T. Filippov, 
A. N. Koryukin, V. N. Zhelyabin, K. N. Ponomarev, V. N. Gerasimov, 
and O. N. Smirnov, for giving careful considetrations to my resalts
notwithstanding the vacation time. Thanks also are due to 
Drs. Jaime Torres Keller and Suemi Rodr\'\i guez-Romo for the 
beautiful facilities for my research work in the Center of
Theoretical Research (FES-C UNAM) and to
Prof. Zbigniew Oziewicz for interesting comments on the subject matter.

\end{document}